\documentclass{gtart_h}  

\def\ifplaintex{\expandafter\ifx\csname documentclass\endcsname\relax}

\ifplaintex 
\else
\fi

\expandafter\ifx\csname beginpicture\endcsname\relax
\expandafter\ifx\csname documentclass\endcsname\relax
\input pictex \else
\input prepictex \input pictex \input postpictex \fi\fi

\def\gt{{\mathsurround=0pt\it $\cal G\mskip-2mu$eometry \&\ 
$\cal T\!\!$opology}}        %  journal title in recommended style

\def\gtp{{\mathsurround=0pt\it $\cal G\mskip-2mu$eometry \&\ 
$\cal T\!\!$opology $\cal P\!$ublications}}  % GT publications

\def\lognumber#1{\def\thelognumber{#1}}
\def\volumenumber#1{\def\thevolumenumber{#1}}
\def\papernumber#1{\def\thepapernumber{#1}}
\def\volumeyear#1{\def\thevolumeyear{#1}}

\def\pagenumbers#1#2{\def\startpage{#1}\def\finishpage{#2}}
\def\published#1{\def\publishdate{#1}}
\def\proposed#1{\def\theproposer{#1}}
\def\seconded#1{\def\theseconders{#1}}
\def\received#1{\def\receiveddate{#1}}
\def\revised#1{\def\reviseddate{#1}}
\def\accepted#1{\def\accepteddate{#1}}

\def\asciiaddress#1{\def\theasciiaddress{#1}}
\def\asciiemail#1{\def\theasciiemail{#1}}

\long\def\asciiabstract#1{\long\def\theasciiabstract{#1}}
\def\asciikeywords#1{\def\theasciikeywords{#1}}

\def\shorttitle#1{\def\theshorttitle{#1}}

\let\\\par\let\thelognumber\relax
\let\thevolumenumber\relax\let\thepapernumber\relax
\let\thevolumeyear\relax\let\thesamplenumber\relax\let\startpage\relax
\let\finishpage\relax\let\publishdate\relax\let\receiveddate\relax
\let\reviseddate\relax\let\accepteddate\relax\let\theasciititle\relax
\let\theasciiauthors\relax\let\theasciiaddress\relax
\let\theasciiabstract\relax\let\theasciikeywords\relax
\let\theasciiemail\relax\let\theshortauthors\relax\let\theshorttitle\relax

\long\def\maketitlep{   % start of definition of \maketitlep

\count0=\startpage

\gt\hfill      %   Journal title (top left) 
\beginpicture
\setcoordinatesystem units <0.33truein, 0.33truein> point at 2.2 0.9
\setplotsymbol ({$\cal G$})
\plotsymbolspacing=9truept
\circulararc 315 degrees from 0 1 center at 0 0
\setplotsymbol ({$\cal T$})
\circulararc 315 degrees from 1 -1 center at 1 0
\endpicture
\break
{\small\ifx\thesamplenumber\relax % sample?  
Volume \else Sample
\fi\thevolumenumber\ (\thevolumeyear)
\startpage--\finishpage\nl
Published: \publishdate}
\vglue 0.5truein plus 0.4fil minus 0.1truein

{\parskip=0pt\leftskip 0pt plus 1fil\def\\{\par\smallskip}{\ifplaintex\large
\else\Large\fi\bf\thetitle}\par\medskip}   

\vglue 0pt plus 0.1fil 

{\parskip=0pt\leftskip 0pt plus 1fil\def\\{\par}{\sc\theauthors}
\par\medskip}

\vglue 0pt plus 0.1fil 

{\small\parskip=0pt\let\newline\\
{\leftskip 0pt plus 1fil\def\\{\par}{\sl\theaddress}\par}
\expandafter\ifx\theemail\relax    % email address?
\relax\else\vglue 5pt plus 0.02fil minus 2pt\def\\{\stdspace{\rm 
and}\stdspace} 
\cl{Email:\stdspace\tt\theemail}\fi
\ifx\theurl\relax                  % URL given?
\relax\else\vglue 5pt plus 0.02fil minus 2pt\def\\{\stdspace{\rm 
and}\stdspace}
\cl{URL:\stdspace\tt\theurl}\fi\par}

\vglue 7pt plus 0.3fil minus 3pt

{\bf Abstract}
\vglue 5pt plus 0.1fil minus 2pt

\theabstract

\vglue 7pt plus 0.3fil minus 3pt

{\bf AMS Classification numbers}\quad Primary:\quad \theprimaryclass

Secondary:\quad \thesecondaryclass

\vglue 5pt plus 0.3fil minus 2pt

{\bf Keywords:}\quad \thekeywords

\vglue 10pt plus 0.5fil minus 5pt

{\small  Proposed: \theproposer\hfill Received: \receiveddate\nl
Seconded: \theseconders\hfill 
\ifx\reviseddate\relax                         % paper revised?
Accepted: \accepteddate                        % no
\else
Revised: \reviseddate                          % yes
\fi}
\eject
}       %  end of definition of \maketitlep

\let\maketitlepage\maketitlep
\let\maketitle\maketitlepage

\font\phead=cmsl9 scaled 950
\font=cmsl9 scaled 1050
\font\pnum=cmbx10 scaled 913
\font\lnum=cmbx10 
\font\pfoot=cmsl9 scaled 950
\font=cmsl9 scaled 1050
\ifplaintex
\headline{\vbox to 0pt{\vskip -4.5mm\line{\small\phead\ifnum
\count0=\startpage ISSN 1364-0380 (on line)
1465-3060 (printed) \hfill {\pnum\folio}\else\ifodd\count0\def\\{ }% 
\ifx\theshorttitle\relax\thetitle\else\theshorttitle\fi\hfill{\pnum\folio}
\else\def\\{ and }{\pnum\folio}\hfill\ifx\theshortauthors\relax\theauthors
\else\theshortauthors\fi\fi\fi}\vss}}
\footline{\vbox to 0pt{\vglue 0mm\line{\small\pfoot\ifnum\count0=\startpage
\copyright\ \gtp\hfill\else
\gt, Volume \thevolumenumber\ (\thevolumeyear)\hfill\fi}\vss
}}
\else
\makeatletter
\def\@oddhead{{\small\lhead\ifnum\count0=\startpage ISSN 1364-0380 (on line)
1465-3060 (printed) \hfill {\lnum\number\count0}\else\ifodd\count0
\def\\{ }\ifx\theshorttitle\relax \thetitle \else\theshorttitle\fi\hfill
{\lnum\number\count0}\else\def\\{ and }{\lnum\number\count0}
\hfill\ifx\theshortauthors\relax 
\theauthors\else\theshortauthors\fi\fi\fi}}\def\@evenhead{\@oddhead}
\def\@oddfoot{\small\lfoot\ifnum\count0=\startpage\copyright\ \gtp\hfill\else
\gt, Volume \thevolumenumber\ (\thevolumeyear)\hfill\fi}
\def\@evenfoot{\@oddfoot}
\makeatother
\fi

\newwrite\gtoutfile
\long\gdef\makeheadfile{  %%% start of definition of \makeheadfile
{\def\\{, }\def\s{ }
\immediate\openout\gtoutfile head.xxx
\immediate\write\gtoutfile{Proxy-for: \ifx\theasciiauthors\relax
\theauthors\else\theasciiauthors\fi\s<\ifx\theasciiemail\relax\theemail\else\theasciiemail\fi>}
\immediate\write\gtoutfile{\noexpand\\}
\immediate\write\gtoutfile{Authors: \ifx\theasciiauthors\relax
\theauthors\else\theasciiauthors\fi}
{\def\\{ }\immediate\write\gtoutfile{Title: \ifx\theasciititle\relax
\thetitle\else\theasciititle\fi}}
\immediate\write\gtoutfile{Subj-class: GT or SG or MG etc}
\immediate\write\gtoutfile{MSC-class: \theprimaryclass\ifx\thesecondaryclass\relax\else, \thesecondaryclass\fi}
\immediate\write\gtoutfile{Journal-ref: Geom. Topol. \thevolumenumber
(\thevolumeyear) \startpage-\finishpage}
\immediate\write\gtoutfile{Comments: Published by Geometry and Topology at}
\immediate\write\gtoutfile{\s\s http://www.maths.warwick.ac.uk/gt/GTVol\thevolumenumber/paper\thepapernumber.abs.html}
\immediate\write\gtoutfile{\noexpand\\}
\immediate\write\gtoutfile{}
\ifx\theasciiabstract\relax
\immediate\write\gtoutfile{\theabstract}\else
\immediate\write\gtoutfile{\theasciiabstract}\fi
\immediate\write\gtoutfile{}
\immediate\write\gtoutfile{\noexpand\\}
\immediate\write\gtoutfile{}
\immediate\closeout\gtoutfile}}  %%% end of definition of \makeheadfile

\def\maketitlepage{\maketitlep\makeheadfile}
\let\maketitle\maketitlepage

\lognumber{408}
\volumenumber{8}\papernumber{20}\volumeyear{2004}
\pagenumbers{743}{777}
\received{17 January 2004}
\accepted{16 May 2004}
\revised{6 May 2004}
\published{18 May 2004}
\proposed{Yasha Eliashberg}
\seconded{Leonid Polterovich, Simon Donaldson}

\usepackage{amsfonts}
\usepackage{amsmath}
\usepackage{graphicx}

\newtheorem{theorem}{Theorem}
\newtheorem{lemma}[theorem]{Lemma}
\newtheorem{proposition}[theorem]{Proposition}
\newtheorem{corollary}[theorem]{Corollary}

\newtheorem{addendum}[theorem]{Addendum}

\theoremstyle{definition}

\newtheorem{example}[theorem]{Example}

\newtheorem{remark}[theorem]{Remark}

\numberwithin{equation}{section}

\def\Z{\mathbb Z}

\def\R{\mathbb R}
\def\C{\mathbb C}
\def\a{\alpha}

\def\e{\epsilon}
\def\l{\lambda}
\def\o{\omega}
\def\s{\sigma}

\def\w{\wedge}
\def\p{\partial}

\newcommand{\into}{\ensuremath{\hookrightarrow}}
\newcommand{\tb}{\mathop{\rm tb}\nolimits}
\newcommand{\rot}{\mathop{\rm rot}\nolimits}
\newcommand{\lk}{\mathop{\rm lk}\nolimits}

\newcommand{\interior}{\mathop{\rm int}\nolimits}

\newcommand{\Span}{\mathop{\rm span}\nolimits}
\newcommand{\genus}{\mathop{\rm genus}\nolimits}

\begin{document}
%%%%

\title{Constructing symplectic forms on 4--manifolds\\which vanish on circles}
\shorttitle{Constructing symplectic forms which vanish on circles}

\authors{David T Gay\\Robion Kirby}
\address{CIRGET, Universit\'{e} du Qu\'{e}bec \`{a} Montr\'{e}al, Case Postale
8888\\Succursale centre-ville, Montr\'{e}al (QC) H3C 3P8, Canada}
\secondaddress{Department of Mathematics, University of California\\Berkeley, CA 94720, USA} 
\asciiaddress{CIRGET, Universite du Quebec a Montreal, Case Postale
8888\\Succursale centre-ville, Montreal (QC) H3C 3P8, 
Canada\\and\\Department of Mathematics, University of 
California\\Berkeley, CA 94720, USA} 
\gtemail{\mailto{gay@math.uqam.ca}{\rm\qua and\qua}\mailto{kirby@math.berkeley.edu}}
\asciiemail{gay@math.uqam.ca, kirby@math.berkeley.edu}

\begin{abstract}
  Given a smooth, closed, oriented 4--manifold $X$ and $\alpha \in
  H_2(X,\Z)$ such that $\alpha\cdot\alpha > 0$, a closed $2$--form $\o$
  is constructed, Poincar\'{e} dual to $\alpha$, which is symplectic
  on the complement of a finite set of unknotted circles $Z$. The
  number of circles, counted with sign, is given by $d= (c_1(s)^2
  -3\sigma(X) -2\chi(X))/4$, where $s$ is a certain $\text{spin}^{\C}$
  structure naturally associated to $\o$.
\end{abstract}
\asciiabstract{% 
Given a smooth, closed, oriented 4-manifold X and alpha in H_2(X,Z)
such that alpha.alpha > 0, a closed 2-form w is constructed,
Poincare dual to alpha, which is symplectic on the complement of a
finite set of unknotted circles. The number of circles, counted with
sign, is given by d = (c_1(s)^2 -3sigma(X) -2chi(X))/4, where s
is a certain spin^C structure naturally associated to w.}

\asciikeywords{Symplectic, 4-manifold, spin^C, almost complex,
  harmonic}

\maketitlepage

\setcounter{section}{-1}
\section{Introduction}

Let $X^4$ be a connected, closed, smooth, oriented 4--manifold with
$b_2^+ >0$. For a Riemannian metric $g$ on $X$, let $\Lambda^2_+$ be
the $3$--plane bundle of self-dual $2$--forms on $X$.  Harmonic
$2$--forms are closed sections of $\Lambda^2_+$; it is
known~\cite{honda} that there exist metrics on $X$ for which there are
harmonic $2$--forms which are transverse to the $0$--section of
$\Lambda^2_+$, so that the $0$--locus is $1$--dimensional and the
$2$--forms are symplectic in the complement of some circles.  Here, we
explicitly construct such $2$--forms and metrics.

Let us say that a connected subset $C$ in a $4$--manifold ``uses up all
the $3$--handles'' if the complement of a regular neighborhood of $C$
has a handlebody decomposition with only $0$--, $1$-- and $2$--handles.

Given $\alpha \in H_2(X;\Z)$ with $\alpha \cdot \alpha > 0$, let
$\Sigma$ be a smoothly imbedded surface in $X$ which uses up all the
$3$--handles and represents $\alpha$.  In Section~\ref{E} we show that such
surfaces exist, and that one can arbitrarily increase the genus of
$\Sigma$ by adding homologically trivial tori while still using up all
the $3$--handles.  Let $c \in H^2(X;\Z)$ be a cohomology class
satisfying  $c \cdot \Sigma = 2-2(\genus(\Sigma) +1) + \Sigma \cdot
\Sigma$;  such $c$ exist because $\genus(\Sigma)$ can be increased if
necessary.  Then choose a $\text{spin}^{\C}$ structure $s$ with
$c_1(s) = c$.

\begin{theorem}[] \label{mainthm}
  Given $\Sigma$ and $s$ as above, there exist a closed $2$--form $\o$
  on $X$, a finite set of signed ($\pm 1$) circles $Z \subset X -
  \Sigma$ (we will show how this sign is natural) bounding disjointly
  imbedded disks, an $\o$--compatible almost complex structure $J$
  on $X - Z$, and a Riemannian metric $g$ on $X$, which satisfy:

\begin{enumerate}

\item $[\o] \in H^2(X;\R)$ is Poincar\'{e} dual to $[\Sigma]$ (and hence
integral).

\item $\o \wedge \o > 0$ on $X-Z$.

\item $\o$ vanishes identically on $Z$.
  
\item There exists a $J$--holomorphic curve $\Sigma'$ which is the
  connected sum of $\Sigma$ with a standard torus in a neighborhood of
  a point in $\Sigma$.

\item There is at least one circle with sign $-1$.
  
\item The sign ($\pm 1$) associated to each circle $Z_i \subset Z$ is
  the obstruction $o(Z_i) \in \pi_3(S^2)=\Z$ to extending $J$ across a
  $4$--ball neighborhood of a certain disk $D_i$ bounded by $Z_i$ and
  the total number of circles in $Z$, counted with sign, is $d=
  ((c_1(s))^2 -3\sigma(X) -2\chi(X))/4$.
  
\item The $\text{spin}^\C$ structure determined by $J|(X-D)$, where
$D$ is the union of the disks $D_i$, is $s$.
  
\item $\o$ is $g$--self-dual (and thus harmonic) and transverse to the
  zero section of $\Lambda^2_+$ (with zero locus $Z$).

\end{enumerate}

\end{theorem}

\begin{remark} \label{a-cx}
The obstruction invariant $o(Z_i)$ depends only on the relative
homology class of the disk $D_i$ in $H_2(X,Z;\Z)$, not on the specific
disk. We will show (proposition~\ref{oZcount}) that $o(Z_i)$ can also
be computed by counting anti-complex (or complex, depending on
orientation) points on $D_i$.

Note that a $\text{spin}^\C$ structure $s$ is an almost complex
structure on the $2$--skeleton of $X$ which extends over the
$3$--skeleton of $X$ ( see page~48 in~\cite{gompf2}).  Then $J$
restricted to the complement of the disks $D_i$ is a $\text{spin}^\C$
structure $s$, and conversely, $J$ may be thought of as an extension
of $s$ across the $2$--disks transverse to the $D_i$.
\end{remark}

Taubes~\cite{taubes1} has initiated a program to study the behavior of
$J$--holomorphic curves in $4$--manifolds equipped with symplectic forms
which vanish along circles, in the hope that this will reveal smooth
invariants of nonsymplectic $4$--manifolds. We in turn hope that our
construction will produce a rich class of examples in which to pursue
this program. In~\cite{scott}, some explicit constructions in terms of
handlebodies are discussed, but not with this much generality and, in
particular, not on closed manifolds. A canonical example that is in
some sense diametrically opposed to our construction is the case of
$S^1 \times Y^3$ for a $3$--manifold $Y$: One chooses an $S^1$--valued
Morse function $f$ on $Y$ with only critical points of index $1$ and
$2$, and defines $\o$ to be $dt \wedge df + \star_3 df$, where
$\star_3$ is the Hodge star operator on $Y$ and $t$ is the $S^1$
coordinate on $S^1 \times Y$. The zero circles are then $S^1 \times p$
for critical points $p \in Y$. Here the zero circles are all
homologically nontrivial; in our construction the zero circles all
bound disks.

\begin{remark} \label{Hondamodel}
  Honda~\cite{honda2} has shown that, given a metric $g$, when a
  harmonic $2$--form $\o$ is transverse to the zero section of
  $\Lambda^2_+$, the behavior of $\o$ on a neighborhood $S^1 \times
  B^3$ near a component of the zero locus $Z$ is given by one of two
  local models. The ``orientable'' model is $\o = dt \w dh + \star_3
  dh$, as in the previous paragraph, where $h\co  B^3 \to \R$ is a
  standard Morse function with a single critical point of index $1$
  (or $2$) at $0$. Thus there is a natural splitting of the normal
  bundle to $Z$ into a $1$--dimensional and a $2$--dimensional bundle.
  The ``nonorientable'' model is a quotient of the orientable model by
  a $Z_2$ action so that the $1$--dimensional bundle becomes a Moebius
  strip. A positive feature of our construction is that the
  nonorientable model never arises.
\end{remark}

\begin{remark}
  In our construction the adjunction inequality is always violated in
  the sense that we have a $J$--holomorphic curve $\Sigma'$ which does
  not minimize genus. When $c_1(s)^2 - 2\chi(X) - 3\s(X) = 0$, we can
  actually cancel the circles in $Z$ at the level of the almost
  complex structure, to get an almost complex structure on all of $X$
  with respect to which a surface which does not minimize genus is
  $J$--holomorphic. This supplements examples of
  Mikhalkin~\cite{mikhalkin} and Bohr~\cite{bohr}. Furthermore, this
  almost complex structure is compatible with a symplectic form
  outside a ball.
  
  If $X$ does in fact support a symplectic structure, at first glance
  it appears that our construction has no hope of recovering that
  fact, since we always produce singular circles and we always
  violate the adjunction inequality. However, there is a different
  type of ``cancellation of singular circles'' that might appear. In
  the $S^1 \times Y^3$ model described above, a flow line connecting
  an index--$1$ critical point which cancels an index--$2$ critical
  point becomes a symplectic cylinder connecting two singular circles
  which can be cancelled symplectically. In our construction, we could
  search for symplectic cylinders with this local model connecting two
  of our circles. To be able to cancel all the circles in this way, at
  least one of these symplectic cylinders would have to intersect
  $\Sigma'$, so that after the cancellation $\Sigma'$ is no longer
  $J$--holomorphic. Furthermore, the change would have to change
  $c_1(s)$, so that $c_1(s) \cdot [\Sigma]$ no longer predicts a minimal
  genus in $[\Sigma]$ of $\genus(\Sigma) + 1$. An interesting project
  is to search for explicit examples of this kind of cancellation.
\end{remark}

\begin{example} \label{3cp2s}
  Consider $X = \#^3 \C P^2$ (which cannot support a symplectic
  structure), with standard generators $\a_1, \a_2, \a_3 \in
  H_2(X;\Z)$ such that $\a_i \cdot \a_j = \delta_{ij}$. Let $s$ be the
  $\text{spin}^\C$ structure for which $c_1(s) \cdot \a_1 = 1$,
  $c_1(s) \cdot \a_2 = 3$ and $c_1(s) \cdot \a_3 = 3$ ($s$ is unique
  because there is no $2$--torsion here). Let $\Sigma$ be
  the standard $\C P^1$ representing $\a_1$ and check that $c_1(s)
  \cdot \a_1 = 2 - 2(\genus(\Sigma)+1) + \a_1 \cdot a_1$.
  
  In our construction we first build $E$, a neighborhood of $\Sigma$,
  as a neighborhood $E'$ of a $J$--holomorphic {\em torus} $\Sigma'$
  together with two extra $2$--handles. This will have a negative
  overtwisted contact structure on its $S^3$ boundary.  Then we build
  $N = X - \interior(E)$ with the standard two $2$--handles each with
  framing $+1$; to attach these two $2$--handles to the $0$--handle
  along Legendrian knots, we need the boundary of the $0$--handle to be
  convex and overtwisted.  One circle in $Z$ is introduced precisely
  to change the standard tight contact structure on the boundary of
  the standard symplectic $0$--handle by a Lutz twist along a
  transverse unknot with self-linking number $-1$ to achieve
  overtwistedness; the second circle is introduced to cancel the first
  circle at the level of almost complex structures, since
  $d=(19-9-10)/4=0$ for this choice of $s$. The second circle
  corresponds to a Lutz twist along a transverse unknot with
  self-linking number $+1$; the self-linking numbers are exactly the
  signs of the circles. If we put the two circles in a single
  $4$--ball, the obstruction to extending $J$ across this ball is zero.
  
  Precisely because the total obstruction to extending $J$ is $d$ ($0$
  in this case), the contact structures on $S^3$ coming from $N$ and
  from $E$ are homotopic. They are also overtwisted and therefore
  isotopic so we can glue $E$ to $N$ symplectically. This finishes
  Example~\ref{3cp2s}.
\end{example}

We can generalize Theorem~\ref{mainthm} to make configurations of
embedded surfaces $J$--holomorphic. Let $\Sigma_1, \ldots, \Sigma_k$ be
smoothly imbedded surfaces in $X$ with pairwise intersections
transverse and positive (self-intersections not necessarily positive)
such that, for each $i$, $\Sigma_i \cdot \Sigma_1+ \ldots + \Sigma_i
\cdot \Sigma_k > 0$. Let $Q$ be the intersection form for a
neighborhood of $\Sigma_1 \cup \ldots \cup \Sigma_k$ and assume that
$\det(Q) \neq 0$. Suppose that $\Sigma_1 \cup \ldots \cup \Sigma_k$
uses up all the $3$--handles. (Again, at the cost of increasing genus
we can use all the $3$--handles.) Let $s$ be a $\text{spin}^\C$
structure on $X$ such that:
\begin{enumerate}
\item $c_1(s) \cdot \Sigma_1 = 2 - 2(\genus(\Sigma_1) + 1) + \Sigma_1
  \cdot \Sigma_1$ and
\item for each $i > 1$, $c_1(s) \cdot \Sigma_i = 2-2\genus(\Sigma_i) +
  \Sigma_i \cdot \Sigma_i$.
\end{enumerate}

\begin{addendum} \label{moresurfaces}
Given $\Sigma_1, \ldots, \Sigma_k$ and $s$ as above, there exist
$\o$, $Z$, $J$ and $g$ as in Theorem~\ref{mainthm} with the following
adjustments to the properties listed in Theorem~\ref{mainthm}:
\begin{enumerate}
\item[\rm(1)] $[\o]$ is Poincar\'{e} dual to $[\Sigma_1] + \ldots +
  [\Sigma_k]$.
\item[\rm(4)] There exists a $J$--holomorphic curve $\Sigma'_1$ which is the
  connected sum of $\Sigma_1$ with a standard torus in a neighborhood
  of a point in $\Sigma_1$.
\item[\rm(4${}'$)] $\Sigma_2, \ldots, \Sigma_k$ are all $J$--holomorphic.
\end{enumerate}

\end{addendum}

\begin{example}
  Let $Y$ be a closed, oriented $3$--manifold with $b_2(Y) > 0$ and let
  $X = S^1 \times Y$. (If $Y$ does not fiber over $S^1$ it is not known, in general, whether $X$ supports a symplectic structure or not.) Let $\Sigma_1$ be a homologically nontrivial
  surface in $Y$ and let $\Sigma_2 = S^1 \times \gamma$, where
  $\gamma$ is a knot in $Y$ transversely intersecting $\Sigma_1$ at
  one point. Choose $s$ such that $c_1(s)$ is Poincar\'{e} dual to $(2
  - 2(\genus(\Sigma_1) + 1)) [\Sigma_2]$. Again we have $d = 0$, so
  our construction gives an almost complex structure which extends
  across a ball containing the two circles in $Z$, and is compatible
  with a symplectic form outside that ball. Also, $\Sigma'_1$ is a
  $J$--holomorphic curve which does not minimize genus.
\end{example}

Addendum~\ref{moresurfaces} also allows us to carry out our
construction on a $4$--manifold which looks like $\R^4$ outside a
compact set, to get a standard symplectic form at infinity. To make
this more precise, let $W$ be a compact, oriented $4$--manifold with
$\p W = S^3$. Let $\Sigma \subset W$ be a properly imbedded surface
which uses up all the $3$--handles, with $\p \Sigma$ unknotted in
$S^3$. Let $[\Sigma]$ refer to the absolute class obtained by capping
off $\p \Sigma$ with a disk in $S^3$, and suppose that 
$[\Sigma] \cdot [\Sigma] > 0$. 
Furthermore suppose there exists an integral lift
$c$ of $w_2(W)$ such that $c \cdot [\Sigma] = -2(\genus(\Sigma)+1) +
[\Sigma] \cdot [\Sigma]$. 
At the cost of increasing $\genus(\Sigma)$, we can always find such a
$c$. 

\begin{corollary}
  There exist a closed $2$--form $\o$ on $W$ which is symplectic on the
  complement of a collection of circles $Z \subset \interior(W)$ and
  $0$ along $Z$, an $\o$--compatible almost complex structure $J$ on
  $W-Z$, and a Riemannian metric $g$ on $W$ with respect to which $\o$
  is self-dual and transverse to $0$, such that:
\begin{enumerate}
\item $[\o]$ is Poincar\'{e} dual to $[\Sigma]$.
\item $\Sigma'$ is $J$--holomorphic, where $\Sigma'$ is the connected
  sum of $\Sigma$ with a standard trivial torus in a ball.
\item $\o| \p W = d\a_0$ for the standard contact form $\a_0$ on $S^3
  = \p B^4 \subset (\R^4,dx_1 \wedge dy_1 + dx_2 \wedge dy_2)$
\end{enumerate}
\end{corollary}

\begin{proof}
Let
\begin{enumerate}
\item $A = (S^2 \times S^2)- B^4$, 
\item $F$ be a properly embedded disk in $A$ normal to one of the $S^2$'s, 
\item $X = W \cup_{S^3} A$, 
\item $\Sigma_1 = \Sigma \cup F$, a closed surface with $\Sigma_1
  \cdot \Sigma_1 = [\Sigma] \cdot [\Sigma]$,
\item $\Sigma_2 = S^2 \times p \subset A$, and
\item $\Sigma_3 = p \times S^2 \subset A$.  
\end{enumerate}
This gives the suitable input for Addendum~\ref{moresurfaces},
with the determinant of the intersection matrix equal to $-[\Sigma]
\cdot [\Sigma]$ and with $s$ chosen so that $c_1(s) \cdot \Sigma_2 =
c_1(s) \cdot \Sigma_3 = 2$. Then the output is standard near $\Sigma_2
\cup \Sigma_3$, so we can restrict to $W$ to get standard behavior
along $\p W$.
\end{proof}

Finally, there are certain special situations where we do not need to
increase the genus of $\Sigma_1$. Let $\Sigma_1, \ldots, \Sigma_k$ be
as in Addendum~\ref{moresurfaces} above, and let $s$ be chosen so
that, for $i =1, \ldots, k$, $c_1(s) \cdot \Sigma_i = 2 -2g(\Sigma_i)
+ \Sigma_i \cdot \Sigma_i$.
\begin{addendum} \label{noextragenus}
Suppose that we are in one of the following three cases:
\begin{enumerate}
\item $k = 2$, both surfaces are spheres, $\Sigma_1 \cdot \Sigma_1
  \geq 2$ and $\Sigma_2 \cdot \Sigma_2 \geq 1$.
\item $k=2$, $\genus(\Sigma_1) =0$, $\genus(\Sigma_2) \geq 1$ and
  $\Sigma_i \cdot \Sigma_i \geq 1$ for $i=1,2$.
\item $k > 2$, $\genus(\Sigma_1) = 0$, $\Sigma_1 \cdot \Sigma_1 \geq
  1$, $\Sigma_1 \cdot \Sigma_2 = 1$ and $\Sigma_1 \cdot \Sigma_i =
  0$ for $i \geq 2$.
\end{enumerate}
Then we have the same conclusions as in Addendum~\ref{moresurfaces}
except that $\Sigma_1$ is $J$--holomorphic rather than $\Sigma'_1$.
\end{addendum}

\rk{Acknowledgements} The authors would like to thank Yasha Eliashberg
for numerous helpful suggestions and questions, Andr\'{a}s Stipsicz
for help understanding homotopy classes of plane fields and Noah
Goodman for some good ideas about making contact structures
overtwisted. We would also like to thank the National Science
Foundation for financial support under FRG grant number
DMS-0244558; the first author's work was also supported by a CRM/ISM
fellowship.

\section{The construction modulo details}

Henceforth all manifolds will be oriented, all symplectic structures
on $4$--mani\-folds will agree with the orientations of the manifolds,
and all contact structures will be co-oriented. However we will work
with both positive and negative contact structures.

Throughout this paper, if we refer to a triple $(\o, J, \xi)$
(appropriately decorated with subscripts) on a pair $(W,\p W)$, we
mean that $W$ is a $4$--manifold with boundary, $\o$ is a closed
$2$--form on $W$ which vanishes along a (possibly empty) collection of
circles $Z$ and is symplectic on $W-Z$, $J$ is an almost complex
structure on $W-Z$, and $\xi$ is a contact structure on $\p W$. For
such triples, $J$ will always be compatible with $\o$, while $\xi$
will be compatible with $\o$ and $J$ in the following sense: There
should exist a Liouville vector field $V$ defined on a neighborhood of
$\p W$ and transverse to $\p W$, inducing a contact form $\a =
\imath_V \o | \p W$, with $\xi = \ker \a$, $J(\xi) = \xi$ and $J(V) =
R_\a$ (the Reeb vector field for $\a$). Note that when $V$ points out
of $W$, $\xi$ is positive (and $\p W$ is said to be convex), and when
$V$ points in, $\xi$ is negative (and $\p W$ is said to be concave).

In this section we will lay out the construction and defer various
details to later sections. First we give the construction starting
with a single surface.

\begin{proof} [The construction (proof of Theorem~\ref{mainthm})] 
    
  Split $X$ as $X = E \cup N$, where $E$ is a $B^2$--bundle
  neighborhood of $\Sigma$ and $N = X - \interior(E)$. Let $Y = \p N =
  - \p E$, let $c = c_1(s)$ and let $g = \genus(\Sigma)$.

Here is a very brief sketch of the construction of $\o$, $Z$ and $J$
on $X$ such that $c_1(J) = c$. We will construct triples $(\o_E, J_E,
\xi_E)$ on $(E,\p E)$ and $(\o_N, J_N, \xi_N)$ on $(N, \p N)$ with the
following properties:
\begin{enumerate}
\item $\o_E$ and $J_E$ are defined on all of $E$.
\item $\o_N$ is defined on all of $N$ but vanishes along $Z \subset
  N$, while $J_N$ is defined on $N - Z$.
\item $\xi_E$ is negative (as a contact structure on $\p E$, hence
  positive on $Y$) and overtwisted.
\item $\xi_N$ is positive and overtwisted.
\item $\xi_E$ is homotopic as a plane field to $\xi_N$.
\item $c_1(J_N) = c|N$ and $c_1(J_E) = c|E$.
\end{enumerate}
Since $\xi_E$ and $\xi_N$ are both overtwisted, positive as contact
structures on $Y$, and homotopic, we know~\cite{eliashberg2} that they
are isotopic as contact structures. This means that, after rescaling
$\o_E$ and perturbing $\o_E$ and $J_E$ on a small collar neighborhood
of $\p E$, we can glue $(\o_E, J_E)$ to $(\o_N, J_N)$ to get $(\o,J)$
on $X$. A little algebraic topology shows that $c_1(J) = c$.

Now we present the construction in more detail.

Choose integers $l_1, \ldots, l_n \in \{-1,+1\}$, with $l_1 = -1$, so
that $\Sigma_{i=1}^n l_i = d$; these will be the signs associated to
the zero circles $Z_1, \ldots, Z_n$. (If $d < 0$ a natural choice is
$n=d$ and $l_1 = \ldots = l_d=-1$ and if $d > 0$ a natural choice is
$n=d+2$, $l_1=-1$ and $l_2 = \ldots = l_{d+2}=1$. However, any choice
will work.)

First we establish a Morse function on $X$ with particular properties.
Consider the obvious Morse function on $E$ with one $0$--handle, $2 g$
$1$--handles and one $2$--handle. Extend this to a Morse function on all
of $X$ which has only $2$--handles, $3$--handles and a single $4$--handle
in $N$, and then turn this Morse function upside down. Introduce a
cancelling $1$--$2$--handle pair inside $N$. Label the $2$--handles
$H^2_1, \ldots, H^2_p,H^2_{p+1}$, where $H^2_{p+1}$ is the single
$2$--handle in $E$ and $H^2_1$ is the $2$--handle from the cancelling
$1$--$2$ pair. Let $H^1_1$ be the $1$--handle cancelled by $H^2_1$.
Slide $H^2_{p+1}$ over $H^2_1$ so that $H^2_{p+1}$ runs over $H^1_1$
once and then (possibly) over some other $1$--handles. Slide $H^1_1$
over these other $1$--handles so as to arrange that $H^2_{p+1}$ runs
only over $H^1_1$.  This gives a new Morse function which we will call
$f$. Note that $f$ still respects the splitting $X = E \cup_Y N$ in
the sense that we may take $Y$ to be a level set. (This is because we
did not slide handles from $N$ over handles from $E$.)

For any given $t \in \R$, let $X_t = f^{-1}(-\infty, t]$ and let $Y_t
= f^{-1}(t) = \p X_t$. For future convenience, reparametrize $f$ so
that $X_k$ is the union of the $0$-- through $k$--handles, with
$N=X_{1.5}$ and $Y = Y_{1.5}$. Also arrange that all the $1$--handles
are in fact in $X_{0.9} \subset X_1$. (We will put the zero locus $Z$
into $f^{-1}[0.9,1]$.) Figure~\ref{F:Morse} illustrates this Morse
function. Note that $f$ gives a cell decomposition to $X$, with the
descending manifolds for each index $k$ critical point being the
$k$--cells. Label the $k$--cells $C^k_i$ (corresponding to handles
$H^k_i$). We will work with cellular cohomology and homology with
respect to this cell decomposition, and represent cohomology classes
by cellular cocycles. Since $X_k$ deformation retracts onto the
$k$--skeleton of $X$, we will frequently represent classes in
$H^i(X_k;\Z)$ by cocycles on the $k$--skeleton of $X$.  Because of the
handle slides in the previous paragraph, we know that $\p C^2_{p+1} =
C^1_1$ in the cellular chain complex coming from $f$.

\begin{figure}[ht!]
\begin{center}
\small
\begin{picture}(0,0)%
\includegraphics{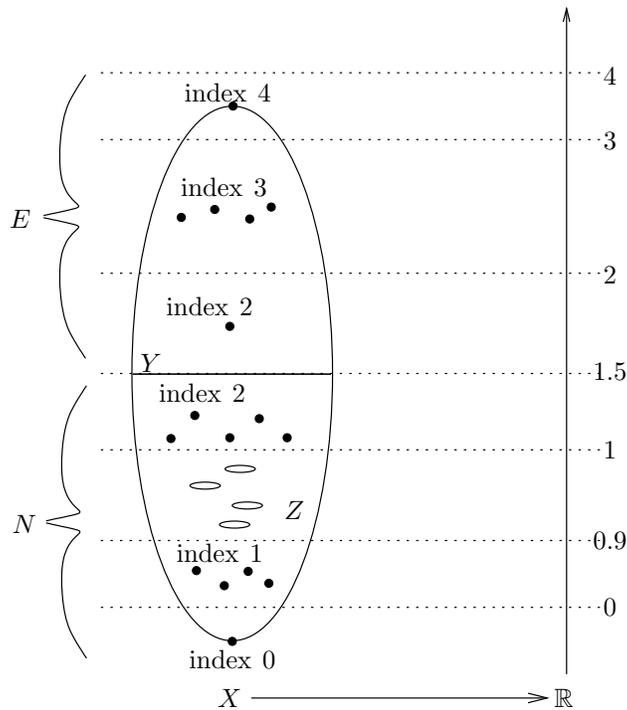}%
\end{picture}%
\setlength{\unitlength}{3947sp}%
\small
\begin{picture}(3670,4437)(2586,-5986)
\put(6000,-5960){\makebox(0,0)[b]{$\R$}%
}
\put(3900,-5960){\makebox(0,0)[b]{$X$}%
}
\put(6290,-5381){\makebox(0,0)[b]{$0$}%
}
\put(6290,-4400){\makebox(0,0)[b]{$1$}%
}
\put(6290,-2046){\makebox(0,0)[b]{$4$}%
}
\put(6290,-2451){\makebox(0,0)[b]{$3$}%
}
\put(6290,-3296){\makebox(0,0)[b]{$2$}%
}
\put(6290,-3906){\makebox(0,0)[b]{$1.5$}%
}
\put(6290,-4966){\makebox(0,0)[b]{$0.9$}%
}
\put(3926,-5715){\makebox(0,0)[b]{index $0$}%
}
\put(3846,-5050){\makebox(0,0)[b]{index $1$}%
}
\put(3736,-4040){\makebox(0,0)[b]{index $2$}%
}
\put(3781,-3500){\makebox(0,0)[b]{index $2$}%
}
\put(3871,-2745){\makebox(0,0)[b]{index $3$}%
}
\put(3896,-2155){\makebox(0,0)[b]{index $4$}%
}
\put(4306,-4775){\makebox(0,0)[b]{$Z$}%
}
\put(2611,-4865){\makebox(0,0)[b]{$N$}%
}
\put(2586,-2945){\makebox(0,0)[b]{$E$}%
}
\put(3406,-3855){\makebox(0,0)[b]{$Y$}%
}
\end{picture}
\caption{The Morse function $f\co X \rightarrow \R$; dots are critical
  points, circles in $f^{-1}[0.9,1]$ are components of $Z$.}
\label{F:Morse}
\end{center}
\end{figure}

In Section~\ref{E} (Proposition~\ref{oEJE}), we construct the triple
$(\o_E,J_E,\xi_E)$. This is done by seeing $E$ as a neighborhood $E'$
of $\Sigma'$ together with two $2$--handles.  The triple is constructed
first on $E'$ as a standard symplectic neighborhood of a
$J$--holomorphic curve $\Sigma'$, with negative but tight contact
boundary, and then the triple is extended over the two extra
$2$--handles required to make $E$ so that the boundary becomes
overtwisted.  Here we use a characterization~\cite{gay2} of
neighborhoods of symplectic surfaces in terms of open books.  Since
$\Sigma'$ is $J_E$--holomorphic we have $c_1(J_E) = c|E$. In
Section~\ref{E} there are various perturbations of the Morse function
inside $E$ (handle slides and introduction of cancelling pairs); after
Section~\ref{E} we will abandon these perturbations and return to the
original Morse function $f$.

Since $X_{0.9}$ is built from a $0$--handle and some $1$--handles, there
is a more or less canonical construction of a triple $(\o_{0.9},
J_{0.9}, \xi_{0.9})$ on $(X_{0.9},Y_{0.9})$, with $\o_{0.9}$ and
$J_{0.9}$ defined everywhere and with $\xi_{0.9}$ positive and tight
(see Proposition~\ref{weinstein}). In Section~\ref{lutz}
(Proposition~\ref{alllutzes}) we show how to extend this to a singular
triple $(\o_1, J_1, \xi_1)$ on $(X_1, Y_1)$, with $\o_1$ vanishing on
a union of circles $Z \subset f^{-1}[0.9,1]$ and $J_1$ defined on
$X_1-Z$, so that $\xi_1$ is positive and overtwisted. In fact, $Z$
consists of one circle in each of $n$ levels between $Y_{0.9}$ and
$Y_1$, in the following sense: Each component $Z_i$ of $Z$ arises in
the product cobordism $f^{-1}[a(i),a(i+1)]$ from $(Y_{a(i)},
\xi_{a(i)})$ to $(Y_{a(i+1)}, \xi_{a(i+1)})$ ($i$ ranging from $1$ to
$n$ and $a(i)$ ranging from $0.9$ to $1$), where each $\xi_{a(i)}$ is
a positive contact structure and $\xi_{a(i+1)}$ differs from
$\xi_{a(i)}$ by a (half) Lutz twist along a transverse unknot $U_i
\subset (Y_{a(i)}, \xi_{a(i)})$. The circle $Z_i$ is $0.5 \times U_i$
after identifying $f^{-1}[a(i),a(i+1)]$ with $[0,1] \times Y_{a(i)}$.
Also, $\o_1|f^{-1}[a(i),a(i+1)]$ is a standard symplectification of
$\xi_{a(i)}$ outside $[0,1] \times T_i$ for some solid torus
neighborhood $T_i$ of $U_i$. Here we also construct a metric $g$ such
that $\o_1$ is $g$--self-dual and transverse to the zero section of
$\Lambda^2_+$.  (Outside a small neighborhood of $Z$, $g$ is
determined by $\o_1$ and $J_1$ in the usual way, but the metric
determined by $\o_1$ and $J_1$ develops a singularity along $Z$ which
we remove by suitably rescaling.)

A homologically trivial transverse knot $K$ comes with an integer
invariant, the self-linking number $\lk(K)$. We will use the fact that
any positive contact manifold has transverse unknots with $\lk = -1$,
and that if the contact structure is overtwisted we can also find
transverse unknots with $\lk = +1$ (lemma~\ref{lk11}). Each $U_i$ is
chosen so that $\lk(U_i) = l_i$ (hence the requirement that $l_1 =
-1$). In Section~\ref{extending} (lemma~\ref{lkequalsobst}) we show
that, if $B_i$ is a $4$--ball neighborhood of a disk bounded by $Z_i$,
then the obstruction to extending $J_1$ across $B_i$ is exactly
$\lk(U_i) = l_i$.

At this point fix a trivialization $\tau$ of $\xi_1$ (possible because
$c_1(\xi_1)=0$). Let $L = K_1 \cup \ldots \cup K_p \cup K_{p+1}$ be
the link of attaching circles for all the $2$--handles $H^2_1 \cup
\ldots \cup H^2_p \cup H^2_{p+1}$ of $X$ as seen in $(Y_1,\xi_1)$.
Each handle $H_i$ is to be attached with some framing $F_i$ of $K_i$.
Because $\xi_1$ is overtwisted, we may isotope $L$ to be Legendrian
with $\tb(K_i)-1 = F_i$ (see lemma~\ref{tbrot}). Thus
(Proposition~\ref{weinstein}, see~\cite{weinstein}), we can extend
$(\o_1, J_1, \xi_1)$ to a triple on $(X_2, Y_2)$ which we will label
$(\o'_2, J'_2, \xi'_2)$. (We will soon change some choices and replace
this triple with a better one, which we will call $(\o_2, J_2,
\xi_2)$.)

With respect to the trivialization $\tau$, each (oriented) Legendrian
knot $K_i$ has a rotation number $\rot(K_i)$ (the winding number of $T
K_i$ in $\xi_1|K_i$ relative to $\tau$). Let $x'$ be the cochain whose
value on a $2$--cell $C^2_i$ is exactly $\rot(K_i)$ for the
corresponding $K_i$. As a cochain on the $2$--skeleton of $X$, $x'$ is
trivially a cocycle; in Section~\ref{cocycles} (lemma~\ref{cocyclec1})
we show that $x'$ represents $c_1(J'_2) \in H^2(X_2;\Z)$ (since $X_2$
deformation retracts onto the $2$--skeleton of $X$). We would like to
have $c_1(J'_2) = c|X_2$, but this is probably not the case ($x'$ is
probably not even a cocycle on $X$).

Represent $c$ by a cocycle $x$ on $X$ which is congruent mod $2$ to
$x'$. (For any representative $x$ of $c$, since both $c_1(J'_2)$ and
$c|X_2$ reduce mod $2$ to $w_2(X_2)$, $x - x'$ is congruent mod $2$ to
$\delta y$ for some $1$--cochain $y$. Thus we can replace $x$ with $x -
\delta y$ if necessary.) In Section~\ref{background}
(lemma~\ref{tbrot}) we show that we can isotope any Legendrian knot
$K$ in an overtwisted contact structure to a new Legendrian knot,
without changing $\tb(K)$, so as to change $\rot(K)$ by any even
number. Thus we can arrange that $\rot(K_i) = x(C^2_i)$ for each
$2$--cell $C^2_i$. Now discard $(\o'_2, J'_2, \xi'_2)$ and use the new
$K_i$'s as attaching circles.  Attach $H^2_1, \ldots, H^2_p$ along
$K_1, \ldots, K_p$ to get $(\o_N, J_N, \xi_N)$ on $(N,Y) =
(X_{1.5},Y_{1.5})$. Then attach $H^2_{p+1}$ along $K_{p+1}$ to extend
$(\o_N,J_N,\xi_N)$ to $(\o_2,J_2,\xi_2)$ on $(X_2,Y_2)$. Thus
(Proposition~\ref{cocyclec1}) $c_1(J_2) = [x] = c|X_2$ and $c_1(J_N) =
[x]|N = c|N$. In the end we will only use $(\o_N,J_N,\xi_N)$, but we
need $J_2$ to show that $\xi_N$ is homotopic to $\xi_E$. Here we use
Proposition~\ref{weinstein} again to know that our construction
extends over the $2$--handles; this proposition also tells us that
$\o_N$ is exact.

To be sure that $\xi_N$ is overtwisted we should arrange that there is
some fixed overtwisted disk which is missed by all the attaching
circles for the $2$--handles. This is possible because, as a result of
a Lutz twist, we have a circles' worth of overtwisted disks, whereas
to adjust $\tb$ and $\rot$ for the attaching circles we only need a
neighborhood of a single overtwisted disk.

We show in Section~\ref{cocycles} (lemma~\ref{cocycleextends}) that,
precisely because $x$ is a cocycle on all of $X$, the almost complex
structure $J_2$ will extend across the $3$--cells of $X$, and hence
will extend to $J_3$ on $X_3$. (Here we abandon the symplectic form
and contact structure.) We still have $c_1(J_3) = c$. Recall the
$4$--ball neighborhoods $B_i$ of disks bounded by the circles $Z_i$.
Let $B_* = f^{-1}[3,\infty)$ be the single $4$--handle for $X$. We can
think of $J_3$ as defined on $X - (B_1 \cup \ldots \cup B_n \cup
B_*)$. It is known~\cite{hh} (and explained in
Section~\ref{extending}, lemma~\ref{totalobstruction}) that the total
obstruction to extending an almost complex structure $J$ defined on
the complement of some balls in a closed $4$--manifold $X$ over those
balls is precisely $(c_1(J)^2 - 3\sigma(X) - 2\chi(X))/4$. But we have
arranged that the sum of the obstructions to extending $J_3$ over
$B_1, \ldots, B_n$ is already $d = (c_1(J_3)^2 - 3\sigma(X) -
2\chi(X))/4$, and thus the obstruction to extending across $B_*$ is
$0$. Therefore $J_3$ extends over the $4$--handle to an almost complex
structure $J_*$ on $X - Z$.

Now compare $J_*|E$ and $J_E$. We know that $c_1(J_*|E) = c|E =
c_1(J_E)$, that $H^2(E;\Z)$ has no $2$--torsion, and that $E$
deformation retracts onto a $2$--complex (coming from the dual Morse
function $-f$). In Section~\ref{planes} (lemma~\ref{Jhomotopy}), we
show that this implies that $J_E$ is homotopic to $J_*|E$. Therefore
$\xi_E$ is homotopic to $\xi_N$ (this also follows from~\cite{gompf})
and we can glue $(E,\o_E, J_E)$ to $(N,\o_N,J_N)$ as described above
to get $(X,\o,J)$.

Now we have $J$ on $X-Z$ with $c_1(J|N) = c|N$ and $c_1(J_E) = c|E$.
This implies that $c_1(J) = c$ because, in the cohomology
Mayer-Vietoris sequence for $X = E \cup N$, the map $H^1(E) \oplus
H^1(N) \rightarrow H^1(E \cap N)$ is surjective, so that the map
$H^2(X) \rightarrow H^2(E) \oplus H^2(N)$ is injective. (The precise
topology of $E$ and $Y = E \cap N$ is important here.)

Now suppose that our construction produced the $\text{spin}^{\C}$
structure $s_0$ when instead we wanted $s_1=s$. We know that $s_1$ is
the result of acting on $s_0$ by some class $a \in H^2(X;\Z)$ of order
$2$.  Our construction above was based on a choice of cocycle
representative $x$ for $c$. In Section~\ref{planes}
(Proposition~\ref{spinCcocycle}), we show that if, instead of $x_0=x$,
we had used $x_1 = x + 2z$ for a special cocycle representative $z$ of
$a$, then we would have produced the desired $\text{spin}^{\C}$
structure $s_1$.

The metric $g$ is constructed first on $X_1$ as mentioned earlier;
on the rest of $X$, $g$ is determined by $\o$ and $J$. In dimension
$4$, if a metric $g$ is given by a symplectic form $\o$ and a
compatible almost complex structure, then $\o$ is automatically
$g$--self-dual.

Finally we can rescale $\o$ so that $\int_{\Sigma} \o = [\Sigma]
\cdot [\Sigma]$. Let $\sigma$ be the Poincar\'{e} dual of $[\Sigma]$
in $H^2(X;\R)$. Then we know that $[\o]|E = \sigma|E$. Also, because
$\o$ is exact on $N$, we know that $[\o]|N = \sigma|N$. As with
$c_1(J)$, the Mayer-Vietoris sequence then gives that $[\o] = \sigma$.
\end{proof}

\begin{proof}[Proof of Addendum~\ref{moresurfaces}]
  
  Now let $E$ be a regular neighborhood of $\Sigma_1 \cup \ldots \cup
  \Sigma_n$ (a plumbing of $B^2$--bundles given by a plumbing graph
  $G$, vertices corresponding to surfaces and edges corresponding to
  intersections), and let $N = X - \interior(E)$.  Recall that $Q$ is
  the intersection form for $E$. The construction is essentially the
  same except for the following changes:
  
  The Morse function we choose on $X$ begins with a standard Morse
  function on $E$ with one $0$--handle for each edge in $G$, one
  $0$--handle for each vertex, one $1$--handle for each incidence
  between an edge and a vertex, and $2\genus(\Sigma_i)$ $1$--handles
  and one $2$--handle for each surface $\Sigma_i$ (see~\cite{gay1}).
  (Of course we can cancel many $0$--$1$--handle pairs, but then the
  picture becomes less canonical.) Extend to $X$ and turn upside down,
  as before. Now, for each $\Sigma_i$, introduce a cancelling
  $1$--$2$--handle pair and slide some handles so that the attaching
  circle for the $2$--handle coming from $\Sigma_i$ runs exactly once
  over exactly one $1$--handle, as we did for the single $2$--handle
  coming from $E$ in the preceding proof.
  
  Given this handlebody decomposition, the construction of
  $(\o_N,J_N,\xi_N)$ is unchanged. In Section~\ref{E}
  (Proposition~\ref{oEJEmoresurfaces}) we show how to adjust the
  construction of $(\o_E,J_E,\xi_E)$ to handle the case where $E$ is a
  plumbing of many disk bundles, again using a characterization of
  neighborhoods of configurations of symplectic surfaces
  in~\cite{gay1}.
  
  To see that our construction gives the correct $c_1(J)$, again we
  use the fact that the map $H^1(E) \rightarrow H^1(Y)$ is surjective
  for the plumbings that we are dealing with. (Here is where we need
  that $\det(Q) \neq 0$.) To see that we can achieve the correct
  $\text{spin}^\C$ structure even in the presence of $2$--torsion, we
  use the same argument from the main proof, since we have arranged
  that each $2$--handle from $E$ runs once over one $1$--handle.

Proposition~\ref{oEJEmoresurfaces} also arranges that $\int_{\Sigma_i}
\o_E = \Sigma_1 \cdot \Sigma_i + \ldots + \Sigma_k \cdot \Sigma_i$ for
each $\Sigma_i$. Thus, if we let $\sigma$ be the Poincar\'{e} dual of
$[\Sigma_1] + \ldots + [\Sigma_k]$, we see again that $[\o]|E =
\sigma|E$ and $[\o]|N = \sigma|N$, so that $[\o] = \sigma$.
\end{proof}

\begin{proof}[Proof of Addendum~\ref{noextragenus}]
  In Section~5.4 of~\cite{goodman}, Goodman investigates conditions
  under which concave boundaries of configurations of symplectic
  surfaces have overtwisted boundaries. The techniques there show
  precisely that the configurations in Addendum~\ref{noextragenus}
  have symplectic neighborhoods with concave overtwisted boundary. We
  use this structure for $(\o_E,J_E,\xi_E)$ and the rest of the
  construction is as before.
\end{proof}

\section{Brief symplectic and contact prerequisites}

\label{background}

Suppose that $W^4$ is a $4$--manifold with boundary $\p W = M^3$, and
suppose that a triple $(\o, J, \xi)$ has been constructed on $(W,M)$.
Let $K$ be a Legendrian knot in $(M,\xi)$. $K$ comes with a natural
contact framing $\tb(K)$ (standing for ``Thurston-Bennequin''), given
by any vector field in $\xi|K$ orthogonal to $K$. Let $W'$ be the
result of attaching a $2$--handle $H$ along $K$ with framing
$\tb(K)-1$. Alternatively, suppose that $W'$ is the result of
attaching a $1$--handle $H$ to $W$, with no constraints on the attaching
map, but with the assumption that $\xi$ is positive. In either case,
we have the following well known result:

\begin{proposition}
[Weinstein~\cite{weinstein}, Eliashberg~\cite{eliashberg3}] 
\label{weinstein}
The triple $(\o, J, \xi)$ extends to $(\o', J', \xi')$ on
$(W', M' = \p W')$, with $\o'$ symplectic on all of $H$ and $J'$
defined on all of $H$. Furthermore $\xi = \xi'$ outside the
surgery that changes $M$ to $M'$. If $\o$ was exact then $\o'$ is
also exact.
\end{proposition}

As originally presented in~\cite{weinstein} and~\cite{eliashberg3},
this applied only to the case where $\o$ was symplectic everywhere and
$\xi$ was positive. As long as the singularities of $\o$ stay in
$\interior(W)$, extending to singular symplectic forms is
trivial. When $\xi$ is negative and $H$ is a $2$--handle, we should
turn the original model handles of~\cite{weinstein}
and~\cite{eliashberg3} upside down. This is discussed
in~\cite{dinggeiges}.

In~\cite{dinggeiges} it is also shown that, in the $2$--handle case,
the surgery that turns $(M,\xi)$ into $(M',\xi')$ is uniquely
determined up to isotopy fixed outside a neighborhood of $K$ by the
property that it preserves tightness near $K$. This unique contact
surgery is called $\tb-1$ surgery along $K$. There is also a uniquely
determined $\tb+1$ surgery, which is precisely what we see if we
attach a $2$--handle as in the proposition to a negative contact
boundary, but reverse orientations so that we see our contact
structure as positive.

When an oriented Legendrian knot $K$ bounds a surface $F \subset M$,
there is a well-defined integer $\rot(K)$, the rotation number of $K$,
given by trivializing $\xi$ over $F$ and counting the winding number
of $TK$ relative to this trivialization. If $c_1(\xi) = 0$, then we
can pick a global trivialization of $\xi$ which works for all $F$'s.
Then, even for homologically nontrivial knots $K$, we get a
well-defined rotation number relative to this trivialization, which we
again call $\rot(K)$. (Note that $\rot(-K) = - \rot(K)$.)

\begin{lemma}
\label{tbrot}
In the above situation, suppose that $a,b$ are integers with $a+b
\equiv \tb(K)+\rot(K) \mod 2$. If $\xi$ is overtwisted then we may
smoothly isotope $K$ to another Legendrian knot $K'$ such that
$\tb(K') = a$ and $\rot(K') = b$.
\end{lemma}

\begin{figure}[ht!]
\begin{center}
\small
\begin{picture}(0,0)%
\includegraphics{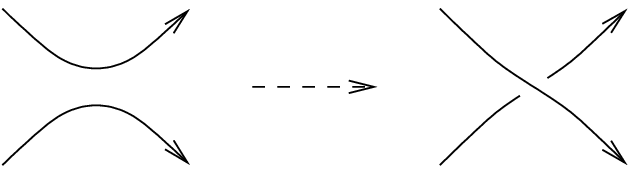}%
\end{picture}%
\setlength{\unitlength}{3947sp}%
\small
\begin{picture}(3024,1041)(1489,-2815)
\put(1951,-2011){\makebox(0,0)[b]{$K_0$}%
}
\put(4051,-2611){\makebox(0,0)[b]{$K_0 \# K_1$}%
}
\put(1951,-2461){\makebox(0,0)[b]{$K_1$}%
}
\end{picture}
\caption{Connected sum of Legendrian knots}
\label{F:legsum}
\end{center}
\end{figure}

\begin{figure}[ht!]
\begin{center}
\includegraphics[width=2in,height=.7in]{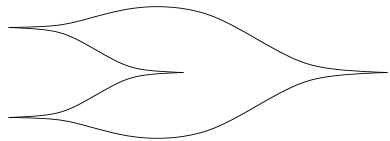}
\caption{Knot with $(\tb,\rot)=(-2,\pm 1)$}
\label{F:tbrotknot}
\end{center}
\end{figure}

\begin{proof}
  In a standard contact chart we can perform a connected sum between
  two Legendrian knots $K_0$ and $K_1$ as illustrated in
  Figure~\ref{F:legsum} (which is a standard ``front projection'' onto
  the $yz$ coordinate plane, where the contact structure is $\ker (dz
  + x dy)$ ). Locally we can assume that our trivialization of $\xi$
  is the vector field $\frac{\p}{\p x}$. Then $\tb(K_0 \# K_1) =
  \tb(K_0)+\tb(K_1)+1$ and $\rot(K_0 \# K_1) = \rot(K_0) + \rot(K_1)$.
  In any positive contact $3$--manifold, there is a Legendrian unknot
  $K_{-2,1}$ with $\tb(K_{-2,1})= -2$ and $\rot(K_{-2,1})=\pm 1$. For
  example, see Figure~\ref{F:tbrotknot}.  In any positive, overtwisted
  contact $3$--manifold there is a Legendrian unknot $K_{0,1}$ with
  $\tb(K_{0,1})=0$ and $\rot(K_{0,1})= \pm 1$. (Consider the contact
  structure $\xi_0 = \ker (\cos(r^2) dz + \sin(r^2) d\theta)$ on
  $\R^3$ with cylindrical coordinates $(r,\theta,z)$, and let $D$ be
  the overtwisted disk $\{r^2 \leq \pi, z = 0\}$.  In every
  overtwisted contact $3$--manifold one can find a ball contactomorphic
  to a neighborhood of $D$ with this contact structure. Let $K_{0,1}$
  be the preimage of $\p D$ under this contactomorphism. One can check
  explicitly that $\tb(\p D) = 0$ and $\rot(\p D) = \pm 1$.)  We
  construct $K'$ as the connected sum of $K$ with some number of
  copies of $K_{-2,1}$ and $K_{0,1}$ so as to arrange that $\tb(K') =
  a$ and $\rot(K')=b$. Of course $K'$ is smoothly isotopic to $K$.
\end{proof}

Now let $K$ be a homologically trivial {\em transverse} knot in
$(M,\xi)$. There is an integer invariant of $K$ known as the
self-linking number $\lk(K)$ given by trivializing $\xi$ over a
surface $F$ bounded by $K$ and seeing this trivialization restricted
to $K$ as a framing of $K$. Again consider the contact structure
$\xi_0 = \ker(\cos(r^2) dz + \sin(r^2) d\theta)$ on $\R^3$. Every
contact $3$--manifold has charts contactomorphic to small neighborhoods
of $0$ in $(\R^3,\xi_0)$, and every overtwisted manifold has charts
contactomorphic to small neighborhoods of $D = \{r^2 \leq \pi\}$.
Consider the circles $K_R = \{r = R\}$.  The following can be easily
checked:
\begin{lemma}
\label{lk11}
For small positive $\e$, $K_\e$ is transverse with $\lk(K_\e) = -1$,
and $K_{\sqrt{\pi}+\e}$ is transverse with $\lk(K_{\sqrt{\pi} + \e}) =
+1$. Thus every contact manifold has knots with self-linking number
$-1$, and every overtwisted contact manifold has knots with
self-linking number $+1$.

\end{lemma}

\section{Constructing $(\o_E, J_E, \xi_E)$}
\label{E}

As a warmup and for the sake of completeness, we prove the following:
\begin{lemma} \label{genus}
  Given any class $\a \in H_2(X;\Z)$, there exists an integer $g(\a)$
  with the following property: For any $g \geq g(\a)$, there is a
  genus $g$ surface $\Sigma$ with $[\Sigma] = \a$ which uses up all
  the $3$--handles.
\end{lemma}

\begin{proof}
  Choose some embedded surface $\Sigma$ representing $\a$.  Let $g =
  \genus(\Sigma)$ and let $W$ be a $B^2$--neighborhood of $\Sigma$.
  There exists a handlebody decomposition of $X$ with one $0$--handle,
  $2g$ $1$--handles and one $2$--handle, forming $W$, together with $l$
  more $1$--handles, some more $2$--handles, some $3$--handles and a
  $4$--handle. Let $H$ be the single $2$--handle in $W$;
  Figure~\ref{F:W} is a picture of $W$. (The framing for $H$ is $\a
  \cdot \a$.)
\begin{figure}[ht!]
\begin{center}
\begin{picture}(0,0)%
\includegraphics{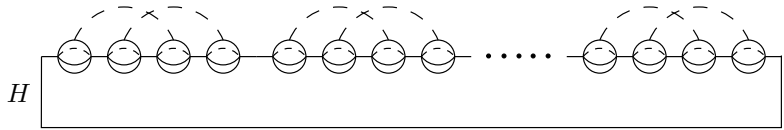}%
\end{picture}%
\setlength{\unitlength}{3947sp}%
\begin{picture}(4700,780)(1463,-2023)
\small
\put(1463,-1844){\makebox(0,0)[rb]{$H$}%
}
\end{picture}
\caption{The initial handlebody decomposition for $W$, with one
  $2$--handle $H$}
\label{F:W}
\end{center}
\end{figure}

Let $g(\a) = g + l$. Given $g' = g(\a)+k \geq g(\a)$, introduce $k$
  cancelling pairs of $1$-- and $2$--handles, so that now we have $q=
  l+k$ extra $1$--handles. 

Let $A_1, \ldots, A_q$ be the extra $1$--handles and introduce $q$ more
cancelling $1$--$2$--handle pairs, with the $1$--handles labelled $B_1,
\ldots, B_q$ and the respective $2$--handles labelled $C_1, \ldots,
C_q$.  Let $W' = W \cup (A_1 \cup \ldots \cup A_q) = W \cup (A_1 \cup
\ldots A_q \cup B_1 \cup \ldots \cup B_q \cup C_1 \ldots \cup C_q)$.
Now, for $i=1$ up to $q$, slide $H$ over $A_i$
then over $C_i$ twice, as in Figure~\ref{F:slides}, to get
Figure~\ref{F:Wprime}. (We have suppressed framings in the figures,
but the $C_i$'s should be $0$--framed so that, after sliding, the
framing on $H$ is still $\a \cdot \a$.)
\begin{figure}[ht!]
\begin{center}
\begin{picture}(0,0)%
\includegraphics{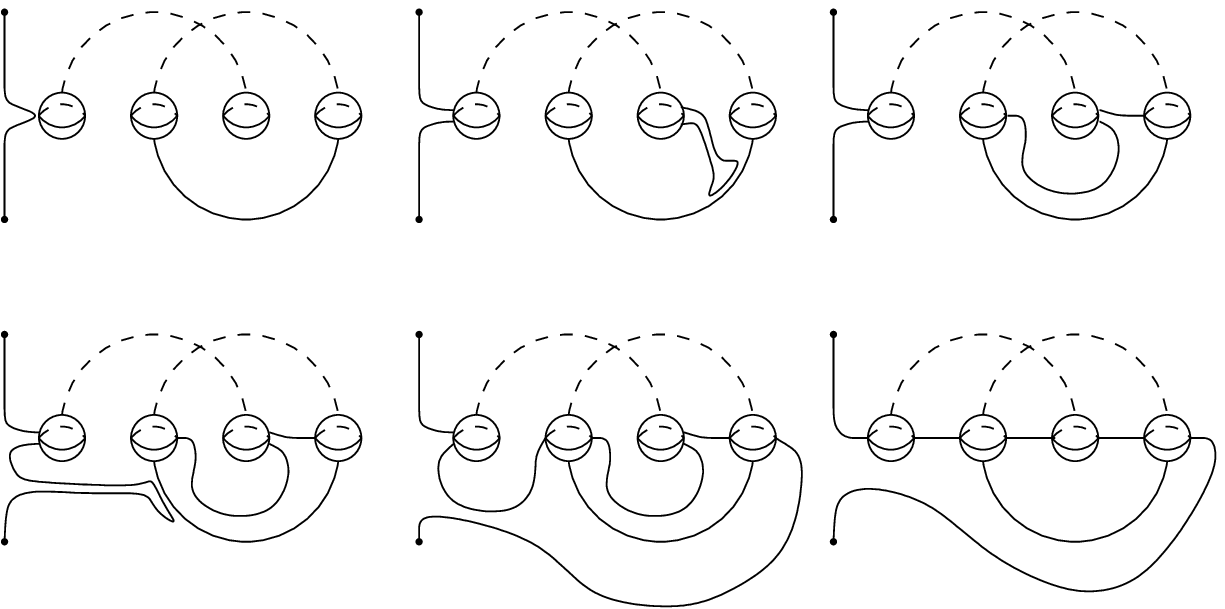}%
\end{picture}%
\setlength{\unitlength}{3947sp}%
\small
\begin{picture}(5845,2992)(160,-2192)
\put(899,699){\makebox(0,0)[b]{$A_i$}%
}
\put(1341,699){\makebox(0,0)[b]{$B_i$}%
}
\put(2889,699){\makebox(0,0)[b]{$A_i$}%
}
\put(3331,699){\makebox(0,0)[b]{$B_i$}%
}
\put(4878,699){\makebox(0,0)[b]{$A_i$}%
}
\put(5320,699){\makebox(0,0)[b]{$B_i$}%
}
\put(3331,-848){\makebox(0,0)[b]{$B_i$}%
}
\put(2889,-848){\makebox(0,0)[b]{$A_i$}%
}
\put(1341,-848){\makebox(0,0)[b]{$B_i$}%
}
\put(899,-848){\makebox(0,0)[b]{$A_i$}%
}
\put(5320,-848){\makebox(0,0)[b]{$B_i$}%
}
\put(4878,-848){\makebox(0,0)[b]{$A_i$}%
}
\put(5320,-488){\makebox(0,0)[b]{$C_i$}%
}
\put(3331,-488){\makebox(0,0)[b]{$C_i$}%
}
\put(1341,-488){\makebox(0,0)[b]{$C_i$}%
}
\put(1341,-2035){\makebox(0,0)[b]{$C_i$}%
}
\put(3331,-2035){\makebox(0,0)[b]{$C_i$}%
}
\put(5320,-2035){\makebox(0,0)[b]{$C_i$}%
}
\put(167,341){\makebox(0,0)[rb]{$H$}%
}
\put(2156,341){\makebox(0,0)[rb]{$H$}%
}
\put(4146,341){\makebox(0,0)[rb]{$H$}%
}
\put(4146,-1151){\makebox(0,0)[rb]{$H$}%
}
\put(2156,-1151){\makebox(0,0)[rb]{$H$}%
}
\put(165,-1151){\makebox(0,0)[rb]{$H$}%
}
\end{picture}
\caption{Sliding $H$ over a $1$--handle then twice over a $2$--handle}
\label{F:slides}
\end{center}
\end{figure}

\begin{figure}[ht!]
\begin{center}
\begin{picture}(0,0)%
\includegraphics{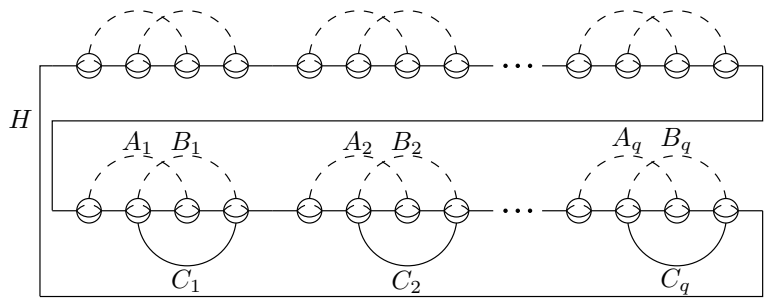}%
\end{picture}%
\setlength{\unitlength}{3947sp}%
\small
\begin{picture}(4580,1840)(-213,-1589)
\put(432,-651){\makebox(0,0)[b]{$A_1$}%
}
\put(740,-651){\makebox(0,0)[b]{$B_1$}%
}
\put(1817,-651){\makebox(0,0)[b]{$A_2$}%
}
\put(2124,-651){\makebox(0,0)[b]{$B_2$}%
}
\put(3509,-651){\makebox(0,0)[b]{$A_q$}%
}
\put(3817,-651){\makebox(0,0)[b]{$B_q$}%
}
\put(-213,-490){\makebox(0,0)[rb]{$H$}%
}
\put(3813,-1535){\makebox(0,0)[b]{$C_q$}%
}
\put(2124,-1525){\makebox(0,0)[b]{$C_2$}%
}
\put(740,-1524){\makebox(0,0)[b]{$C_1$}%
}
\end{picture}
\caption{Handlebody decomposition of $W'$}
\label{F:Wprime}
\end{center}
\end{figure}

With this new handlebody decomposition let $W'' = W' - (C_1 \cup
\ldots \cup C_q)$, and note that $W''$ is a neighborhood of a surface
$\Sigma'$ of genus $g + q = g(\a) + k = g'$. The remainder of $X$ is
built with $C_1, \ldots, C_q$ and some more $2$--handles, $3$--handles
and a $4$--handle. Thus $X - W''$ has a handlebody decomposition with
only $0$--, $1$-- and $2$--handles.

To see that $[\Sigma'] = [\Sigma]$, note that we slid $H$ over each
$C_i$ once in a positive direction and once in a negative direction.
\end{proof}

\begin{figure}[ht!]
\begin{center}
\begin{picture}(0,0)%
\includegraphics{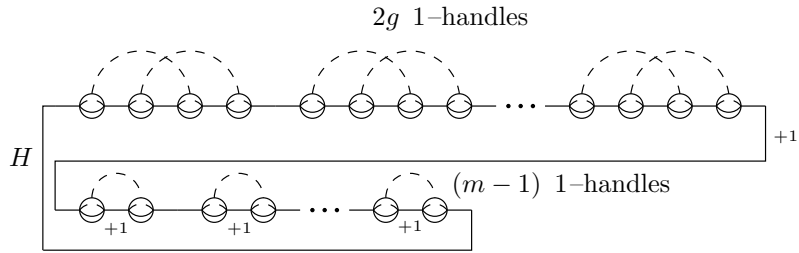}%
\end{picture}%
\setlength{\unitlength}{3947sp}%
\small
\begin{picture}(4624,1542)(-230,-1023)
\put(-230,-485){\makebox(0,0)[rb]{$H$}%
}
\tiny
\put(277,-918){\makebox(0,0)[b]{$+1$}%
}
\put(1051,-918){\makebox(0,0)[b]{$+1$}%
}
\put(2124,-907){\makebox(0,0)[b]{$+1$}%
}
\put(4394,-336){\makebox(0,0)[lb]{$+1$}%
}
\small
\put(2374,-689){\makebox(0,0)[lb]{$(m-1)$ $1$--handles}%
}
\put(1875,375){\makebox(0,0)[lb]{$2g$ $1$--handles}%
}
\end{picture}
\caption{Handlebody decomposition of $E$}
\label{F:E}
\end{center}
\end{figure}

Now we return to the notation in the introduction: $\Sigma$ is a
given surface of genus $g$ representing a class $\a \in H_2(X;\Z)$,
with $B^2$--neighborhood $E$. Let $m = \a \cdot \a > 0$. The Morse
function $-f\co  X \to \R$ restricts to the obvious Morse function on $E$
with one $0$--handle, $2 g$ $1$--handles and one $2$--handle, which we
will label $H$ as in the proof of the preceding lemma. First introduce
$m-1$ cancelling $1$--$2$--handle pairs inside $E$ (each new $2$--handle
framed $+1$) and slide $H$ over
all the new $2$--handles so that $E$ gets a handlebody decomposition
with one $0$--handle, $(2g + m - 1)$ $1$--handles and $m$ $2$--handles
attached as in Figure~\ref{F:E}. Now the framing on $H$ will also be $+1$.

\begin{figure}[ht!]
\begin{center}
\begin{picture}(0,0)%
\includegraphics{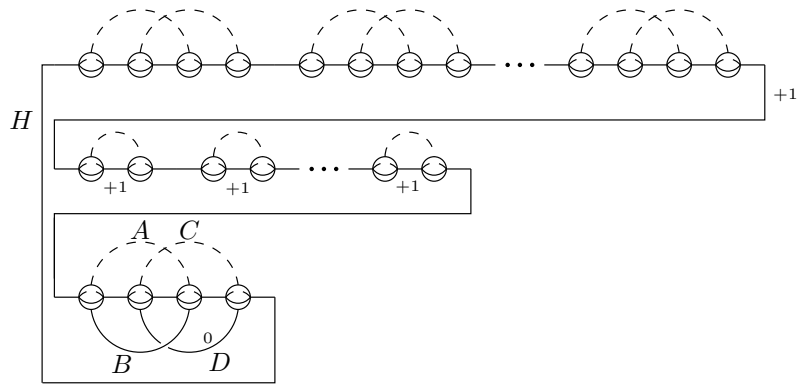}%
\end{picture}%
\setlength{\unitlength}{3947sp}%
\begin{picture}(4618,2366)(-224,-2115)
\tiny
\put(4394,-336){\makebox(0,0)[lb]{$+1$}%
}
\small
\put(432,-1199){\makebox(0,0)[b]{$A$}%
}
\put(740,-1199){\makebox(0,0)[b]{$C$}%
}
\put(930,-2032){\makebox(0,0)[b]{$D$}%
}
\put(317,-2043){\makebox(0,0)[b]{$B$}%
}
\tiny
\put(277,-918){\makebox(0,0)[b]{$+1$}%
}
\put(1045,-924){\makebox(0,0)[b]{$+1$}%
}
\put(2113,-912){\makebox(0,0)[b]{$+1$}%
}
\small
\put(-224,-507){\makebox(0,0)[rb]{$H$}%
}
\tiny
\put(857,-1855){\makebox(0,0)[b]{$0$}%
}
\end{picture}
\caption{Another handlebody decomposition of $E$}
\label{F:E2}
\end{center}
\end{figure}

Next introduce a $1$--handle $A$ cancelled by a $2$--handle $B$ and a
$1$--handle $C$ cancelled by a $2$--handle $D$, inside $E$. As in the
proof of the preceding lemma, as in Figure~\ref{F:slides}, slide $H$
over $A$ and then over $D$ twice, to get Figure~\ref{F:E2}. To get the
right picture, the framing of $D$ should be $0$ but the framing of $B$
can be anything so we have left it unlabelled in the figure.
With respect to this new handlebody decomposition, let $E' = E - (B
\cup D)$.  Note that $E'$ is a neighborhood of a surface $\Sigma'$
of genus $g' = g+1$ formed as the connected sum of $\Sigma$ with
a standard torus in a $4$--ball neighborhood of a point in $\Sigma$.

The first important point is that, if we attach $D$ to $E'$ then we
can slide $H$ back off of $D$ (twice) and over $A$ and then cancel $C$
and $D$. Therefore $E' \cup D = E \natural (S^1 \times B^3)$, where
the $S^1 \times B^3$ summand comes precisely from $A$.  Thus $\p (E'
\cup D) = \p E \# (S^1 \times S^2)$, and $B$ is to be attached along
any circle isotopic to $S^1 \times \{p\}$ in the $S^1 \times S^2$
summand.

The second important point is an understanding of the framings of the
cancelling $2$--handles $B$ and $D$. When constructing $E'$ there is a
natural open book decomposition of $\#^{(2g'+m-1)} (S^1 \times S^2)$
(the boundary of the $0$--handle and the $1$--handles) with page an
$m$--punctured genus $g'$ surface and trivial monodromy. The
construction of $E'$ is completed by attaching $m$ $2$--handles along
the binding, each with framing $+1$ (relative to the page), so that
$\p E'$ inherits a natural open book decomposition. The page can still
be seen clearly in Figure~\ref{F:E2} as the disk bounded by the
attaching circles for the $2$--handles (not including $B$ and $D$)
together with $2$--dimensional $1$--handles embedded in the
$4$--dimensional $1$--handles. The attaching circles for $B$ and $D$ can
be isotoped to lie on this page. To make the sliding of $H$ work
right, we needed $D$ to be attached with framing $0$ measured relative
to this page.

To put a symplectic structure on $E$ with concave, overtwisted
boundary, we will need the following:

\begin{lemma} \label{L:OTsurgery}
  Given any {\em negative} contact structure $\xi$ on a $3$--manifold
  $M^3$ and any Legendrian knot $K$ which transversely intersects a
  $2$--sphere $S$ in $M$ at one point, $\tb(K)-1$ surgery along $K$
  leads to an overtwisted contact structure.
\end{lemma}

\begin{proof}
  First reverse the orientation of $M$. Then $\xi$ is a positive
  contact structure and we need to show that $\tb(K)+1$ surgery along
  $K$ leads to an overtwisted contact structure. We will show that,
  after surgery, $S$ becomes a disk with Legendrian boundary $K'$ (the
  dual circle to $K$) with $\tb(K')=+1$, which is not possible in a
  positive tight contact structure.
  
  Let $T_0 = D^2 \times S^1$ be a solid torus with coordinates
  $(r,\mu,\lambda)$, where $(r,\mu)$ are polar coordinates on $D^2$
  and $D^2 = \{r \leq 1\}$. Every Legendrian knot has a neighborhood
  which is contactomorphic to $T_0$ with a certain standard contact
  structure $\xi_0$ with the property that $\p T_0 = S^1 \times S^1$
  is a convex surface with dividing set $\Gamma_0 = \{r=1, \mu \in
  \{0,\pi\}\}$ (see~\cite{dinggeiges}, for example). Here $K_0 =
  \{r=0\}$ is Legendrian and $\tb(K_0)$ is the $0$--framing coming from
  our splitting of $T_0$ as $D^2 \times S^1$. To perform surgery along
  a Legendrian knot $K \subset (M,\xi)$, find a neighborhood $T$ of
  $K$ contactomorphic to $(T_0,\xi_0)$. Then $M - \interior(T)$ has
  convex torus boundary with dividing set $\Gamma$ equal to two
  parallel longitudes. Now glue $(T_0,\xi_0)$ back in via any
  diffeomorphism $\phi\co  \p T_0 \rightarrow \p T$ which takes
  $\Gamma_0$ to $\Gamma$. Identifying $\p T$ with $\p T_0$ via the
  original contactomorphism between $(T,\xi)$ and $(T_0,\xi_0)$, we
  think of $\phi$ as an automorphism of $\p T_0$ taking $\Gamma_0$ to
  $\Gamma_0$, and hence we think of $\phi \in SL(2,\mathbb{Z})$.
  Legendrian $\tb+1$ surgery corresponds to $\phi = \left(
  \begin{array}{cc}1 & 0 \\ 1 & 1 \end{array} \right)$, using the
  basis $(\mu, \lambda)$.
  
  Recall our $2$--sphere $S$: Note that $S \cap \p T$ is a meridian,
  or, after identifying $T$ with $T_0$ and lifting to the universal
  cover, the line spanned by $v = (1,0)^T$. Then $\phi^{-1}(v) =
  (1,-1)^T$. In other words, after surgery, $D = S - \interior(T)$ is
  a disk which meets $\p T_0$ in a longitude representing the framing
  $\tb(K_0)-1$. Extend $D$ in to $T_0$ by the obvious annulus so that
  $\p D = K_0 = K'$. Then the (topological) canonical zero-framing of
  $K'$ given by $D$ is $\tb(K')-1$, hence $\tb(K') = +1$.
\end{proof}

\begin{proposition}\label{oEJE}
There exists a triple $(\o_E, J_E, \xi_E)$ on $(E, \p E)$ such that
$\xi_E$ is negative and overtwisted, $\Sigma'$ is $J_E$--holomorphic
and $c_1(J_E) = c|E$.
\end{proposition}

\begin{proof}
  
  Since $E'$ is a neighborhood of a surface $\Sigma'$ with positive
  self-inter\-section, $E'$ has a symplectic structure $\o_{E'}$ with
  concave boundary (by Theorem~1.1 in~\cite{gay2}), in which $\Sigma'$
  is symplectic. We can find a compatible almost complex structure
  $J_{E'}$ with respect to which $\Sigma'$ is holomorphic, and thus we
  get our triple $(\o_{E'}, J_{E'}, \xi_{E'})$, with $\xi_{E'}$
  negative.  Theorem~1.1 in~\cite{gay2} also says that $\xi_{E'}$ is
  supported by the open book on $\p E'$ described above. (See
  Section~2 of~\cite{gay1} and Section~1 of~\cite{gay2} for background
  on the important relationship between contact structures and open
  books established by Giroux.)  Let $K_D$ be the attaching circle for
  the $2$--handle $D$; we noted earlier that $K_D$ lies in a page of
  the open book.  Since $K_D$ is also homologically nontrivial in that
  page, we may assume that $K_D$ is Legendrian, with $\tb(K_D)$ equal
  to the framing coming from the page, which in this case means
  $\tb(K_D)=0$ (see, for example, Remark~4.1 in~\cite{gay1}, or the
  Legendrian realization principle of~\cite{honda3}).  Since
  $\xi_{E'}$ is a {\em negative} contact structure, we can easily
  isotope $K_D$ to another Legendrian knot so as to {\em increase}
  $\tb(K_D)$, so that $\tb(K_D) = 1$, and hence $\tb(K_D)-1 = 0$, the
  desired framing for $D$. Thus $(\o_{E'}, J_{E'}, \xi_{E'})$ extends
  over $(E' \cup D, \p (E' \cup D))$ (see
  Proposition~\ref{weinstein}).
  
  In $\p (E' \cup D)$, the attaching circle $K_B$ of $B$ transversely
  intersects a $2$--sphere at one point; make $K_B$ Legendrian with
  this property and attach $B$ along this Legendrian knot with framing
  $\tb(K_B)-1$, to get $(\o_E,J_E,\xi_E)$ (again using
  Proposition~\ref{weinstein}). (In the $1$--handle $A = B^1 \times B^3$, the 
  required $2$--sphere is $p \times \partial B^3$ for any $p\in B^1$.) By 
  lemma~\ref{L:OTsurgery}, $\xi_E$ is overtwisted.
  
  Since $\Sigma'$ is symplectic, we may choose an $\o_E$--compatible
  almost complex structure $J_E$ such that $\Sigma'$ is
  $J_E$--holomorphic. Since $E'$ is a neighborhood of $\Sigma'$,
  $c_1(J_E)|E'$ is determined by its action on $[\Sigma']$, which is
  $2-2g' + \alpha \cdot \alpha = 2-2(g+1) + \a \cdot \a$. But
  $i^*\co H^2(E;\Z) \rightarrow H^2(E';\Z)$ is an isomorphism, and thus
  $c_1(J_E)= c|E$.
\end{proof}

Now recall the more general setting in Addendum~\ref{moresurfaces}: We
have a configuration of surfaces $\Sigma_1, \ldots, \Sigma_k$,
intersecting transversely and positively, such
that $\Sigma_1 \cdot \Sigma_i + \ldots + \Sigma_k \cdot \Sigma_i >
0$. Let $E$ be a plumbed neighborhood of $\Sigma_1 \cup \ldots \cup
\Sigma_k$, let $\Sigma'_1$ be the connected sum of $\Sigma_1$ with a trivial
torus and let $E'$ be a plumbed neighborhood of $\Sigma'_1 \cup \ldots
\cup \Sigma_k$. There is a natural open book on $\p E'$ and two
$2$--handles $B$ and $D$ which, when attached to $E'$, give back
$E$. Locally near $B$ and $D$ the picture looks exactly like
Figure~\ref{F:E2}. 

\begin{proposition}\label{oEJEmoresurfaces}
  There exists a triple $(\o_E, J_E, \xi_E)$ on $(E, \p E)$ such that
  $\xi_E$ is negative and overtwisted, $\Sigma'_1, \Sigma_2, \ldots,
  \Sigma_k$ are all $J_E$--holomorphic and $c_1(J_E) = c|E$.
  Furthermore we can arrange that the $\o_E$ area of each $\Sigma_i$
  is $\Sigma_1 \cdot \Sigma_i + \ldots + \Sigma_k \cdot \Sigma_i$.
\end{proposition}

\begin{proof}
  Theorem~1.1 in~\cite{gay2} gives a construction of symplectic forms
  on neighborhoods of configurations exactly of the type we are
  dealing with here. We use this to construct $E'$ with concave
  boundary; Theorem~1.1 of~\cite{gay2} also says that the contact
  structure on $\p E'$ is supported by the natural open book mentioned
  above, so that $B$ and $D$ can be attached to get an overtwisted
  contact boundary for $E$. Furthermore the areas of each $\Sigma_i$
  produced in~\cite{gay2} are precisely the areas stated above, and
  the surfaces $\Sigma'_1, \Sigma_2, \ldots, \Sigma_k$ are symplectic
  and can therefore be made $J$--holomorphic.
\end{proof}

As mentioned in the introduction, we should now abandon all
perturbations of our Morse function that have been introduced in this
section, and return for the rest of this paper to the original Morse
function $f$.

\section{Lutz twist cobordisms}\label{lutz}

In this section we show how to construct a singular symplectic form on
a product cobordism which connects a contact structure $\xi_0$ at the
bottom to a contact structure $\xi_1$ at the top, with $\xi_1$ being
the result of changing $\xi_0$ by a Lutz twist along some transverse
knot.

Let $K$ be a transverse knot in some positive contact $3$--manifold.
Since $K$ is transverse, it has a solid torus neighborhood which is
contactomorphic to the following model. (Any two transverse knots have
contactomorphic neighborhoods due to a standard
Moser-Weinstein-Darboux argument.)

On the solid torus $\e B^2 \times S^1$, choose polar coordinates
$(r,\mu ), 0 \leq r \leq \e$, on $\e B^2$, and $\l$ on $S^1$. It is
convenient to reparameterize by letting $\rho = r^2/2$ so that $d\rho
= rdr$ and $rdrd\mu d\l = d\rho d\mu d\l = dV$, the volume form for
the solid torus.

\begin{figure}[ht!]
\begin{center}
\begin{picture}(0,0)%
\includegraphics{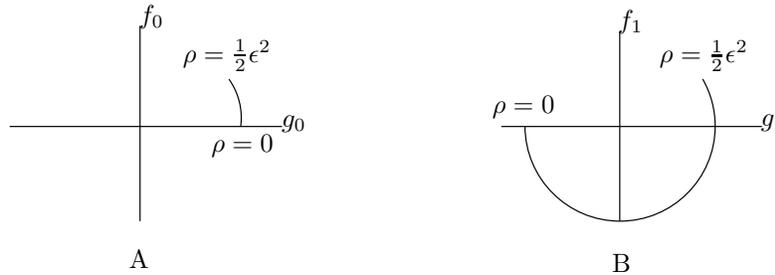}%
\end{picture}%
\setlength{\unitlength}{3947sp}%
\small
\begin{picture}(4809,1651)(1189,-4174)
\put(2987,-3295){\makebox(0,0)[b]{$g_0$}%
}
\put(2014,-4151){\makebox(0,0)[b]{A}%
}
\put(5045,-4174){\makebox(0,0)[b]{B}%
}
\put(5998,-3284){\makebox(0,0)[b]{$g_1$}%
}
\put(2664,-3446){\makebox(0,0)[b]{$\rho = 0$}%
}
\put(2567,-2909){\makebox(0,0)[b]{$\rho = \frac{1}{2}\e^2$}%
}
\put(4429,-3202){\makebox(0,0)[b]{$\rho = 0$}%
}
\put(5561,-2911){\makebox(0,0)[b]{$\rho = \frac{1}{2}\e^2$}%
}
\put(2092,-2643){\makebox(0,0)[b]{$f_0$}%
}
\put(5102,-2672){\makebox(0,0)[b]{$f_1$}%
}
\end{picture}
\caption{Graphs of $\rho \mapsto (g_0(\rho), f_1(\rho))$ and $\rho
  \mapsto (g_1(\rho), f_1(\rho))$}
\label{F:fg}
\end{center}
\end{figure}

Consider a $\mu \l $--invariant form $\a_0 = f_0(\rho ) d\mu + g_0(\rho
) d\l $.  Then $d\a_0 = f_0^{\prime}d\rho d\mu + g_0^{\prime}d\rho
d\l$ and then $\a_0 \wedge d\a_0 = (g_0 f_0^{\prime} - f_0
g_0^{\prime} ) dV > 0$ implies $\frac{d}{d\rho} (\frac{f_0}{g_0}) > 0$.
The contact planes, namely $\xi_0 = \text{ker}\a_0 = \text{span}\{
\frac{\partial}{\partial\rho}, g_0(\rho )\frac{\p}{\p \mu} - f_0(\rho )
\frac{\p}{\p \l} \}$, should be orthogonal to the core circle $0
\times S^1$ parameterized by $\l$, so $f_0(0) = 0$ and $g_0(0) = 1$ is a
good choice, and the map $\rho \to (g_0(\rho), f_0(\rho))$ looks
qualitatively like that drawn in Figure~\ref{F:fg}A.

The graph of $(g_0(\rho), f_0(\rho))$ in Figure~\ref{F:fg}A indicates that
the contact planes, orthogonal to $0\times S^1$, rotate in a left
handed fashion as they move radially away from $0 \times S^1$
(left-handed because $\{ \frac{\partial}{\partial\rho}, g_0(\rho
)\frac{\p}{\p \mu} - f_0(\rho ) \frac{\p}{\p \l} \}$ span the planes).
We can assume they rotate only slightly as $\rho$ traverses $[0,
\frac{1}{2}\e^2 ]$.

Introducing a Lutz twist about $0 \times S^1$ means that we change the
contact structure rel $\p \e B^2 \times S^1$ so that the contact
planes rotate left-handedly an extra $\pi$ as $\rho$ runs through $[0,
\frac{1}{2}\e^2 ]$, starting at $\rho = 0$ orthogonal to $0 \times
S^1$ but with the opposite orientation, and ending up in the same
position as in the standard model above when $\rho = \frac{1}{2}\e^2$
(see Figure~\ref{F:fg}B). This produces a new contact form $\alpha_1 =
f_1(\rho) d\mu + g_1(\rho) d\lambda$, which we may assume exactly
equals $\alpha_0$ for $\rho$ near $\e^2/2$. Note that some authors
would call this a ``half Lutz twist''.

On the trivial bordism $I \times \e B^2 \times S^1$ the standard
symplectization of $\a_0$ is $d(e^t\alpha_0 ) = e^t (dt \wedge
\alpha_0 + d\alpha_0 ), t \in I$.

\begin{proposition}  \label{locallutz} 
On $I \times \e B^2 \times S^1$, there exists a closed 2-form $\o$ and
a metric $g$ satisfying

\begin{enumerate}

\item $\o \equiv 0 $ on $ Z = \frac{1}{2} \times 0 \times S^1$,

\item  $\o \wedge \o > 0$ on the complement of $Z$,

\item $\o = d(e^t \a_0 )$ on a neighborhood of 
$0\times \e B^2 \times S^1$ and $I \times \p \e B^2 \times S^1$,

\item $\o = d(e^t \a_1 )$ on a neighborhood of 
$1 \times \e B^2 \times S^1$,

\item $\o$ is self-dual with respect to $g$ and transverse to the zero
  section of $\Lambda^2_+$.

\item $\o$ and $g$ define an $\o$--compatible almost complex structure outside
  a small neighborhood of $Z$.

\end{enumerate}
\end{proposition}

\begin{figure}[ht!]
\begin{center}
\begin{picture}(0,0)%
\includegraphics{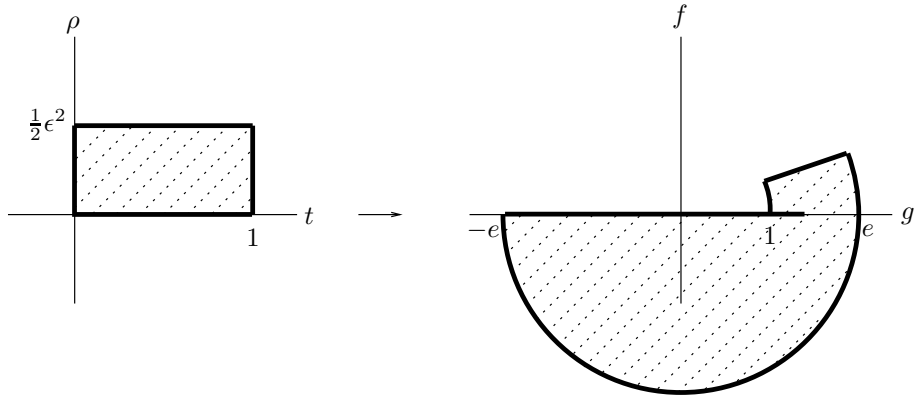}%
\end{picture}%
\setlength{\unitlength}{3947sp}%
\small
\begin{picture}(5669,2466)(289,-3678)
\put(4528,-1347){\makebox(0,0)[b]{$f$}%
}
\put(1845,-2713){\makebox(0,0)[b]{$1$}%
}
\put(5094,-2706){\makebox(0,0)[b]{$1$}%
}
\put(5708,-2652){\makebox(0,0)[b]{$e$}%
}
\put(714,-1377){\makebox(0,0)[b]{$\rho$}%
}
\put(2199,-2576){\makebox(0,0)[b]{$t$}%
}
\put(3291,-2673){\makebox(0,0)[b]{$-e$}%
}
\put(5958,-2577){\makebox(0,0)[b]{$g$}%
}
\put(688,-2050){\makebox(0,0)[rb]{$\frac{1}{2}\e^2$}%
}
\end{picture}
\caption{The map $\phi$}
\label{F:immphi}
\end{center}
\end{figure}

\begin{proof} 

Let $\phi\co  [0,1] \times [0,\e^2/2] \into \R^2$ be a smooth
map satisfying the following properties (see Figure~\ref{F:immphi})
(we use the coordinates $(t,\rho)$ on the domain and $(g,f)$ on the
range for reasons which will become clear shortly):
\begin{enumerate}
\item $\phi$ is an orientation preserving immersion away from
  $(t=1/2, \rho=0)$.
\item $\phi(t,0) = (g(t,0),0)$ for all $t \in [0,1]$.
\item On a neighborhood of $0 \times [0,\e^2/2]$ and
  $[0,1] \times \frac{1}{2} \e^2$, 
  \[\phi(t,\rho) =  e^t(g_0(\rho),f_0(\rho))\] 
  where $g_0$ and $f_0$ are as in the
  preceding paragraphs.
\item On a neighborhood of $1 \times [0,\e^2/2]$, $\phi(t,\rho) =
  e^t(g_1(\rho),f_1(\rho)$.
\item $\phi | [0,\frac{1}{2}) \times 0$ moves monotically from $(1,0)$
  to $(1+\delta,0)$ for some $\delta > 0$, and $\phi |(\frac{1}{2}, 1]
  \times 0$ moves monotonically from $(1 + \delta,0)$ to $( -e, 0)$.
\item In a neighborhood of $(1/2,0) \mapsto (1+\delta,0)$,
  \[ \phi(t,\rho) = (\rho - (t-\frac{1}{2})^2 + 1+\delta, 
                     -2\rho (t-\frac{1}{2})) .\]
\end{enumerate}
(The specified map near $(1/2,0) \mapsto (1+\delta,0)$ behaves much
like the complex map $z \mapsto z^2$, and in particular folds a half-disk
neighborhood of $(1/2,0)$ in the upper half plane onto a disk
neighborhood of $(1+\delta,0)$.)

Letting $t$ be the coordinate on $I$, we have coordinates
$(t,\rho=r^2/2, \mu, \lambda)$ on $I \times \e B^2 \times S^1$.
Writing $\phi(t,\rho) = (g(t,\rho),f(t,\rho))$, let $\a = f(t,\rho)
d\mu + g(t,\rho) d\lambda$, which is a $1$--form on $I \times \e B^2
\times S^1$. (The fact that $f(t,0)=0$ means that $\a$ is well-defined
along $0 \times 0 \times S^1$.) Finally, let $\o = d\a$. The fact that
$\phi$ is an orientation preserving immersion away from $(1/2,0)$
implies that $\o \wedge \o > 0$ away from $Z$. The fact that $d\phi =
0$ at $(1/2,0)$ implies that $\o \equiv 0$ along $Z$. The given
boundary conditions for $\phi$ give the announced boundary conditions
for $\o$. 

Now we must construct $g$ and verify self-duality, transversality, and
compatibility with $\o$. Near $Z$, we have
\[\omega = d\alpha = -2\rho dt
\w d\mu - 2t d\rho \w d\mu - 2t dt \w d\lambda + d\rho \w d\lambda .\]
Now convert polar coordinates $(r,\mu)$ (recalling that $\rho =
r^2/2$) back to cartesian coordinates $(x,y)$ on $\e B^2$ and let $T =
t-1/2$, to get that $\omega = y (dT \w dx + dy \w d\lambda) - x (dT \w
dy - dx \w d\lambda) -2T (dx \w dy + dT \w d\lambda)$. With respect to
the flat metric $g_0 = dT^2 + dx^2 + dy^2 + d\lambda^2$, the three
sections $A= dT \w dx + dy \w d\lambda$, $B = dT \w dy - dx \w
d\lambda$ and $C = dx \w dy + dT \w d\lambda$ give a frame for the
bundle of self-dual $2$--forms, and thus $\omega = y A - x B -2T C$ is
self-dual and transverse to the zero section.

Let $R = \sqrt{4T^2 + x^2 + y^2}$ and let $f(R)$ be a smooth, positive
function which equals $R$ for $R \geq \e'$ and equals $1$ for $R \leq
\e'/2$, for some small $\e'> 0$. Let $g = f(R) g_0$. Note that $\o$ is
still $g$--self-dual and transverse to the zero section of
$\Lambda^2_+$ (since $g=g_0$ near $Z$). But now $g$ and $\o$ induce an
$\o$--compatible almost complex structure $J$ for $R \geq \e'$, by $J =
\tilde{g}^{-1} \circ \tilde{\o}$. (Here $\tilde{g}$ and $\tilde{\o}$
are the maps from the tangent space to the cotangent space induced by
$g$ and $\o$.) In fact, in local coordinates $(T,x,y,\lambda)$ on $R
\geq \e'$, $\tilde{g} = R I$ (where $I$ is the identity matrix) and
$\tilde{\o}$ can be calculated from the explicit form of $\o$ in the
preceding paragraph, to get that $\tilde{\o}$ is the matrix:
\[ Q = \left[ \begin{array}{cccc}
0 & -y & x & 2T \\
y & 0 & 2T & -x \\
-x & -2T & 0 & -y \\
-2T & x & y & 0 
\end{array} \right]
\]
Then $J = (1/R) Q$, $J^2 = (1/R)^2 Q^2 = -I$ and $\omega(Jv,Jw) = v^T
J^T Q J w$, with $J^T Q J = -(1/R)^2 Q^3 = Q$, so that $\o(Jv,Jw) =
\o(v,w)$.
\end{proof}

\begin{remark}
In the proof above, we can rewrite $\omega$ near $Z$ as $d\lambda \w
(2T dT - x dx - y dy) + (-2T dx \w dy - x dT \w dy + y dT \w dx) =
d\lambda \w dh + \star_3 dh$, where $h = -\frac{1}{2} x^2 -
\frac{1}{2} y^2 + T^2$ and $\star_3$ is the Hodge star operator on
$\R^3$ with the flat metric $dy^2 + dx^2 + dT^2$. (Note that $(y,x,T)$
is the correct orientation for $\R^3$ here because we are now writing
our $4$--manifold as $S^1 \times \R^3$.)  This is exactly the oriented
local model given by Honda~\cite{honda2}.
\end{remark}

\begin{remark}
  If $\xi_0$ is a negative contact structure, we can also perform a
  Lutz twist along a transverse knot to get an overtwisted negative
  contact structure $\xi_1$.  One might expect a similar singular
  symplectic cobordism from $\xi_0$ on the bottom to $\xi_1$ on the
  top. Upside down, this would be a cobordism from a positive contact
  structure to a positive contact structure which eliminates a Lutz
  twist. There does exist a singular almost complex structure on such
  a cobordism, but it is not clear how to construct a singular
  symplectic form with the desired properties. Much of this paper
  would be simplified if such a $2$--form could be constructed.
\end{remark}

Now recall the notation from the main construction in the
introduction.
\begin{proposition} \label{alllutzes}
There exists a triple $(\o_1, J_1, \xi_1)$ on $(X_1,
Y_1)$, with the following properties:
\begin{enumerate} 
\item $\o_1$ vanishes on a union of circles $Z \subset
f^{-1}[0.9,1]$.
\item $J_1$ is defined on $X_1 - Z$.
\item $\xi_1$ is positive and overtwisted.
\item $Z$ consists of one unknotted circle $Z_i$ in each
of $n$ levels between $Y_{0.9}$ and $Y_1$. 
\item The obstruction to extending $J_1$ across a $4$--ball
  neighborhood $B_i$ of a disk bounded by $Z_i$ is $l_i$.
\item The metric defined by $\o_1$ and $J_1$ can be modified
  in a small neighborhood of $Z$ to be a metric on all of $X_1$
  such that $\o$ is $g$--self-dual and transverse to the zero section of
$\Lambda^2_+$.
\end{enumerate}
\end{proposition}

\begin{proof}
  First build $(\o_{0.9}, J_{0.9}, \xi_{0.9})$ on $(X_{0.9},Y_{0.9})$
  using standard symplectic $0$-- and $1$--handles, as discussed in
  Section~\ref{background}, Proposition~\ref{weinstein}. Choose
  numbers $0.9 = a(1) < \ldots < a(n+1) = 1$. Choose $K_1 \subset
  (Y_{0.9},\xi_{0.9})$ to be a transverse unknot with $\lk(K_1) = l_1
  = -1$. First put a standard symplectification of $\xi_{0.9}$ on
  $f^{-1}[0.9,a(2)] \cong [0,1] \times Y_{0.9}$, and then replace this
  symplectification with a singular symplectic form on $[0,1] \times
  T_1$, as constructed in Proposition~\ref{locallutz}, for a small
  neighborhood $T_1$ of $K_1$. (This is possible because the
  symplectification and the singular form constructed in
  Proposition~\ref{locallutz} agree on $\p ([0,1] \times T_1)$.) This
  gives $\xi_{a(2)}$ on $Y_{a(2)}$ which is overtwisted. Now we can
  choose $K_2$ a transverse unknot in $(Y_{a(2)},\xi_{a(2)})$ to have
  $\lk(K_2) = l_2$, and repeat. We will prove the statement about the
  obstruction to extending $J_1$ in the next section
  (lemma~\ref{lkequalsobst}). Note that the metrics coming from
  Proposition~\ref{locallutz} fit together smoothly because, away from
  $Z$, they are defined by the symplectic and almost complex
  structures, which fit together smoothly.
\end{proof}

\begin{remark}
There is an alternate construction for $(\o_1, J_1, \xi_1)$ as
follows. Begin with the disjoint union of several copies of $S^1
\times B^3$, on each of which we put the singular symplectic form $dt
\wedge df + \star_3 df$ for $f = -x^2 + y^2 + z^2$ on $B^3$. The
boundary is convex and overtwisted. Connect these together with
$1$--handles, preserving convexity on the boundary, and then kill each
$S^1$ with a Legendrian $2$--handle to get $B^4$, and then attach extra
$1$--handles as necessary. Then one needs to understand how the
obstruction to extending $J_1$ across various balls depends on the
Legendrian $2$--handles used to kill the $S^1$'s. In the end this
construction should be equivalent to that outlined above.
\end{remark}

\section{Extending almost complex structures}
\label{extending}

Consider an almost complex structure $J$ on $TS^3 \oplus \e^1$ where
$\e^1$ is a trivial line bundle which we identify with the normal
bundle $\nu$ to $S^3 = \p B^4$. Also let $\nu$ denote the outward unit
normal vector field spanning the bundle $\nu$. We define four
invariants of $J$ up to homotopy:
\begin{enumerate}
\item Trivialize $TS^3$ using a right-invariant quaternionic frame and thus
  identify unit vectors in $TS^3$ with points in $S^2$. With
  respect to this trivialization, $J(\nu)$ then gives a map $S^3 \to S^2$, and
  thus an element $h(J) \in \pi_3(S^2) = \Z$.
\item Using coordinates $(t,x,y,z)$ on $B^4$, $J(\frac{\p}{\p t})$ gives a map
  $S^3 \to S^2$, where $S^2$ is now the unit $(x,y,z)$ sphere. This gives an
  element $h'(J) \in \pi_3(S^2) = \Z$.
\item Using coordinates $(t,x,y,z)$ on $B^4$, let $u, v, w$ be any
  field of frames for $\Span(\frac{\p}{\p x}, \frac{\p}{\p y},
  \frac{\p}{\p z})$. Now, at each point in $S^3$, interpret
  $J(\frac{\p}{\p t})$ as a point in the $(u,v,w)$ sphere, giving
  another element $h''(J) \in \pi_3(S^2) = \Z$.
\item Choose some $4$--manifold $W$ with $\p W = S^3$ such that $J$
  extends over $W$. This gives the invariant $\theta(J) = (c_1(J)^2 -
  2 \chi(W) - 3 \s(W))/4$.
\end{enumerate}

\begin{remark} \label{additivity}
  The invariant $h'(J)$ is simply the obstruction to extending $J$
  across $B^4$. Thus if we have an almost complex structure defined on
  the complement of two balls $B_1$ and $B_2$ in a $4$--manifold $X$,
  and $B$ is a ball containing $B_1$ and $B_2$, then $h'(J|\p B) =
  h'(J|\p B_1) + h'(J|\p B_2)$. Conversely, if $J$ is an almost
  complex structure on $\p B$ and we choose two integers $k_1, k_2$
  with $k_1 + k_2 = h'(J)$, we can put two balls $B_1$ and $B_2$ in
  $B$ and extend $J$ across $B - (B_1 \cup B_2)$ so that $h'(J|\p B_i)
  = k_i$.
\end{remark} 

\begin{lemma}
These invariants are related by $h(J) = h'(J) = h''(J)$ and $\theta(J)
= -h(J)-\frac{1}{2}$ and each uniquely characterizes $J$ up to
homotopy.
\end{lemma}

\begin{proof}
A direct computation shows that $h(J) = h'(J)$. Since the field of frames
$(\frac{\p}{\p x},\frac{\p}{\p y},\frac{\p}{\p z})$ is homotopic to any other
field of frames $(u,v,w)$ on $B^4$, we get that $h'(J)=h''(J)$.

Note that $J$ defines an oriented plane field $\xi$ on $S^3$, the field of
$J$--complex tangencies to $S^3$. The homotopy class of $J$ is uniquely
determined by the homotopy class of $\xi$.

That $\theta(J) = - h(J) - \frac{1}{2}$ follows from Section~4
of~\cite{gompf}: We are really looking at invariants of $\xi$;
$\theta(J)$ is Gompf's $\theta(\xi)$. The invariant $h(\xi) = h(J)$
can be defined with respect to any trivialization of $TS^3$. Whichever
trivialization we choose, there is a canonical $\Z$ action on homotopy
classes of oriented plane fields which adds $1$ to $h(\xi)$, and Gompf
proves (Theorem~4.5 in~\cite{gompf}) that adding $1$ to $h(\xi)$
corresponds to subtracting $1$ from $\theta(\xi)$. Let $\xi_0$ be the
standard tight positive contact structure on $S^3$. Direct calculation
shows that $h(\xi_0) = 0$ (for our particular trivialization of
$TS^3$) and that $\theta(\xi_0) = -\frac{1}{2}$. Thus $\theta(J) = -
h(J) - \frac{1}{2}$.

Finally it is well known that $h(\xi)$ is a complete invariant for
homotopy classes of oriented plane fields on $S^3$.
\end{proof}

\begin{lemma}
\label{totalobstruction}
Given any closed $X^4$ with balls $B_1, \ldots, B_n \subset X$ and an almost
complex structure $J$ on $X - (B_1 \cup \ldots \cup B_n)$, let $d(J) =
(c_1(J)^2 - 2\chi(X) - 3\s(X))/4$. Then $d(J) = \Sigma_{i=1}^n h(J|\p B_i)$.
\end{lemma}

\begin{proof}
  By Remark~\ref{additivity}, we can assume that $k_i = h(J|\p B_i) =
  \pm 1$ for each $i$. Let $X' = X \#^n \C P^2$, formed by replacing
  each $B_i$ with the complement of a ball in $\C P^2$. Let $E_i$ be
  the $i$'th new generator in $H_2(X';\Z)$ coming from the $i$'th $\C
  P^2 - B^4$. On the $i$'th $\C P^2 - B^4$, put an almost complex
  structure $J_i$ such that $c_1(J_i) \cdot E_i = - 2 k_i + 1$. We
  claim that $J$ and $J_i$ are homotopic on $\p B_i$, so that we can
  glue them together to get an almost complex structure $J'$ on all of
  $X'$.  This is true because $\theta(J_i|\p B_i) = - k_i -\frac{1}{2}
  = - h(J|\p B_i) -\frac{1}{2} = \theta(J|\p B_i)$.
  
  Now, since $J'$ is defined on all of $X'$, we know (by~\cite{hh})
  that $0 = d(J') = d(J) - \Sigma_{i=1}^n k_i = d(J) - \Sigma_{i=1}^n
  h(J|\p B_i)$, and hence $d(J) = \Sigma_{i=1}^n(h|\p B_i)$.
\end{proof}

Now we want to understand the obstruction $o(Z_i)$ to extending $J$
across a ball containing a component $Z_i$ of the singular locus $Z$
for an almost complex structure $J$ coming from our Lutz twist
construction. Suppose $U$ is a transverse unknot in a $3$--manifold
with positive contact structure $\xi$, with a chosen vector field $w$
normal to $\xi$ and tangent to $U$. Let $D$ be a $2$--disk bounded by
$K$ and let $B$ be a $3$--ball neighborhood of $D$. Let $\xi'$ be the
result of performing a Lutz twist along $U$, with vector field $w'$
which is normal to $\xi'$, tangent to $U$, and agrees with $w$ outside
a neighborhood of $U$.  The relevant almost complex structure $J$ is
(up to homotopy) defined on $\p (I \times B)$ as follows: On $(0
\times B) \cup (I \times \p B)$, $J(\frac{\p}{\p t}) = w$ and $J$
preserves $\xi$, and on $1 \times B$, $J(\frac{\p}{\p t}) = w'$ and
$J$ preserves $\xi'$.

\begin{lemma}
\label{lkequalsobst}
With definitions as above, $o(\{0.5\} \times U) = h(J) = \lk(U)$.
\end{lemma}

\begin{proof}
  Recall that $\lk(U)$ is determined by choosing any nonzero section
  $u$ of $\xi|D$ and measuring the framing of $U$ given by $u|U$
  relative to the canonical zero-framing of $U$ coming from $D$.
  
  Extend $u$ to a section of $\xi$ on all of $B$ and let $v = Ju$.  We
  will compute $h''(J)$ using the frame $(u,v,w)$ on $B$ and the
  Thom-Pontrjagin construction. At each point $(t,p)$ on $\p (I \times
  B)$, we write $J_{(t,p)}(\frac{\p}{\p t}) = a(t,p) u + b(t,p) v +
  c(t,p) w$, normalized so that $||(a,b,c)|| = 1$, giving a map $\phi
  \co  (t,p) \mapsto (a(t,p),b(t,p),c(t,p)) \in S^2$.  Let $q = (0,0,-1)$
  and $q' = (-0.1,0,-1)$ normalized to be in $S^2$. We want to compute
  the framed cobordism class of $L = \phi^{-1}(q)$ framed by $L' =
  \phi^{-1}(q')$. On $(0 \times B) \cup (I \times \p B)$, $\phi$ maps
  everything to $(0,0,1)$.  On $1 \times B$, $L$ is exactly $U$, and
  $L'$ is a parallel copy of $U$ realizing the framing given by $u$.
  Thus $h(J) = h''(J) = \lk(U)$.
\end{proof}

More generally, suppose we have any harmonic $2$--form $\o$ on $X$
(with respect to some metric $g$) which is transverse to $0$ and that
$Z$ is a single unknotted component of the zero locus, with the
orientable local model for $\o$ near $Z$ (see
Remark~\ref{Hondamodel}). Now we do not have in mind a particular
contact $3$--manifold $M$ near $Z$ or a Morse function with $Z$ lying
in a regular level set. However we can still choose a $2$--disk $D$
bounded by $Z$ (disjoint from the rest of the zero locus) and a
$4$--ball neighborhood $B$ of $D$ and let $o(Z,D)$ be the obstruction
to extending $J$ across $B$.  Although it is not necessary for this
construction, it would be nice to see how to compute $o(Z,D)$ by
looking directly at the behavior of $J$ along $D$. 

For this we prescribe the manner in which $D$ approaches $Z$, since
$\o$ induces a natural splitting of the normal bundle to $Z$ into $1$--
and $2$--dimensional sub-bundles. Choose local, oriented coordinates
$(\theta,u,v,w)$ near $Z$, with $\theta \in S^1$ and $(u,v,w) \in
B^3$, such that $\o = d\theta \w df + \star_3 df$, $f=
(-u^2-v^2+w^2)/2$ and $\star_3$ is the Hodge star operator on $B^3$
with metric $du^2 + dv^2 + dw^2$ and orientation $(u,v,w)$. Note that
our choice to make $f$ have index $2$ at $0$ fixes an orientation of
$Z$, given by $\theta$. Let $D$ be an oriented, imbedded $2$--disk with
$\p D = Z$ with respect to this orientation, such that, near $Z$, $D$
coincides with the annulus $\{(\theta,u,v,w) | v=w=0,u \geq 0\}$. Note
that $D$ has a singular foliation with singularities corresponding to
complex and anticomplex points (points $p$ where $J(T_p D) = T_p D$
and $J$ agrees with or disagrees with the orientation of $D$). The
foliation can be defined by splitting $TX|D$ as $\nu D \oplus TD$
($\nu D$ being the normal bundle to $D$), choosing a section $V$ of
$\nu D$, homotoping $J$ to be compatible with a product metric on $\nu
D \oplus TD$, and projecting $J(V)$ onto $TD$ to get a vector field on
$D$ which we then integrate to get an oriented singular
foliation. Thus, generically, both complex and anticomplex points can
have either elliptic or hyperbolic neighborhoods. With our boundary
conditions on $D$, we will see that $D$ is complex on a neighborhood
of $\p D$, so that the foliation is not generic there. However, we may
make the foliation generic away from this neighborhood. With respect
to this foliation, let $e_-$ be the number of elliptic anticomplex
points and let $h_-$ be the number of hyperbolic anticomplex points.

\begin{proposition} \label{oZcount}
In the above situation, $o(Z,D) = -1 + 2(e_- - h_-)$.
\end{proposition}

\begin{proof}
A model for a neighborhood of $D$ is $W = D^2_{1+\e} \times D^2_\e$
(where $D^2_r$ is the disk of radius $r$ and $S^1_r = \p D^1_r$), with
$D = D^2_1 \times \{0\}$ and $Z = S^1_1 \times \{0\}$. Let $(x,y)$
be cartesian coordinates on $D^2_{1+\e}$ and $(z,t)$ be cartesion
coordinates on $D^2_\e$. Let $(r,\theta)$ be polar coordinates on
$D^2_{1+\e}$. The coordinates in the discussion above on a
neighborhood of $Z$ are then $(\theta, u=1-r, v=z, w=t)$. Assume $J$
is compatible with the metric $dx^2 + dy^2 + dz^2 + dt^2$ and let
$\phi\co  (W - Z) \rightarrow S^2$ be the map given by seeing $J(\p_t)$
as a point in the unit $(\p_x,\p_y,\p_t)$--sphere.

Consider the following subsets of $W$ (see Figure~\ref{F:subsets}):
\begin{enumerate}
\item $A=\{(r,\theta,z,t) | 1-\e \leq r \leq 1+\e\}$
\item $B=\{(r,\theta,z,t) | 0 \leq r \leq 1-\e, z=t=0\}$
\item $C=\p A \cap \p W$
\item $E = \{(r,\theta,z,t) | 0 \leq r \leq 1-\e, \sqrt{z^2+t^2} =
  \e\} = D^2_{1-\e} \times S^1_\e \subset \p W$
\item $F = \{(r,\theta,z,t) | r=1-\e\} = S^1_{1-\e} \times D^2_\e
  \subset \p A$.
\end{enumerate}

\begin{figure}[ht!]
\begin{center}
\small
\begin{picture}(0,0)%
\includegraphics[width=.9\hsize]{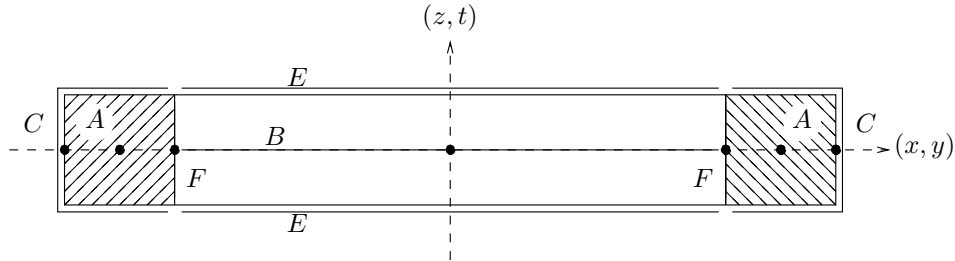}%
\end{picture}%
\setlength{\unitlength}{3670sp}%{3947sp}
\begin{picture}(6022,1650)(1477,-2698)
\small
\put(1737,-1819){\makebox(0,0)[rb]{$C$}%
}
\put(7190,-1819){\makebox(0,0)[lb]{$C$}%
}
\put(3348,-1509){\makebox(0,0)[lb]{$E$}%
}
\put(3348,-2500){\makebox(0,0)[lb]{$E$}%
}
\put(2667,-2190){\makebox(0,0)[lb]{$F$}%
}
\put(6260,-2190){\makebox(0,0)[rb]{$F$}%
}
\put(7450,-2004){\makebox(0,0)[lb]{$(x,y)$}%
}
\put(4464,-1137){\makebox(0,0)[b]{$(z,t)$}%
}
\put(3283,-1906){\makebox(0,0)[b]{$B$}%
}
\put(2067,-1796){\makebox(0,0)[b]{$A$}%
}
\put(6841,-1796){\makebox(0,0)[b]{$A$}%
}
\end{picture}
\caption{Various labelled subsets of $W$}
\label{F:subsets}
\end{center}
\end{figure}

Note that $\p W = C \cup E$. 
Let $\pi\co  \p W \rightarrow B \cup C \cup F = B \cup \p A$ be the map
defined as follows: On $C$, $\pi$ is the identity map. On $E$, $\pi$
maps $\{1-2\e \leq r \leq 1-\e\}$ onto $F$, $\{1-3\e \leq r \leq
1-2\e\}$ onto $\{1-\e \leq r \leq 1-3\e, z=t=0\} \subset B$ and $\{0
\leq r \leq 1-3\e\}$ onto $\{0 \leq r \leq 1-3\e,z=t=0\} \subset B$,
as indicated in Figure~\ref{F:pi}. Let $\psi = \phi \circ \pi\co  \p W
\rightarrow S^2$, and note that $\psi$ is homotopic to $\phi|\p W$. 
\begin{figure}[ht!]
\begin{center}
\small
\begin{picture}(0,0)%
\includegraphics[width=.9\hsize]{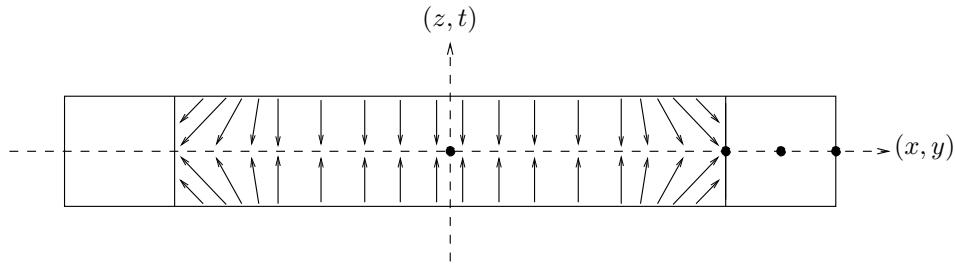}%
\end{picture}%
\setlength{\unitlength}{3670sp}%
\begin{picture}(6022,1650)(1477,-2698)
\small
\put(7450,-2004){\makebox(0,0)[lb]{$(x,y)$}%
}
\put(4464,-1137){\makebox(0,0)[b]{$(z,t)$}%
}
\end{picture}
\caption{The projection $\pi\co  \p W \rightarrow B \cup C \cup F$}
\label{F:pi}
\end{center}
\end{figure}

We will compute $o(Z,D)$ as the oriented framed cobordism class of the
oriented link $\psi^{-1}(0,0,-1)$. From the local form $\o = d\theta
\w df + \star_3 df$ near $Z$, with $f = (-u^2 - v^2 + t^2)/2 =
(-(1-r)^2 - z^2 + t^2)/2$, we compute that, on $A$, $J(\p_t) = (z \p_r
- t\p_\theta + (1-r) \p_z)/\rho$, where $\rho = \sqrt{(1-r)^2 + z^2 +
t^2}$. Thus, on $\p A$, the set of points where $J(\p_t) = -\p_z$ is
just $L_0 = \{(r,\theta,z,t) | r=1+\e, z=t=0\} = S^1_{1+\e} \times
\{0\} \subset C \subset \p W$. Furthermore, as $(\phi|\p
A)^{-1}(0,0,-1) = \psi^{-1}(0,0,-1) \cap C$, $L_0$ has framing $-1$
and is oriented in the negative $\theta$ direction.

\begin{figure}[ht!]
\begin{center}
\small
\begin{picture}(0,0)%
\includegraphics{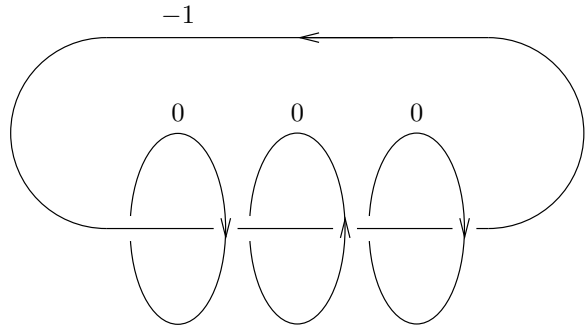}%
\end{picture}%
\setlength{\unitlength}{3947sp}%
\begin{picture}(3618,1989)(592,-3370)
\small
\put(1651,-1486){\makebox(0,0)[b]{$-1$}%
}
\put(1651,-2086){\makebox(0,0)[b]{$0$}%
}
\put(2401,-2086){\makebox(0,0)[b]{$0$}%
}
\put(3151,-2086){\makebox(0,0)[b]{$0$}%
}
\end{picture}
\caption{$\psi^{-1}(0,0,-1)$, in the case of two elliptic anticomplex
  points and one hyperbolic anticomplex point}
\label{F:Pontrjagin}
\end{center}
\end{figure}

On $E$, $\psi^{-1}(0,0,-1)$ is $\pi^{-1}(\phi^{-1}(0,0,-1))$. Since
$\phi^{-1}(0,0,-1) \cap F = \emptyset$, we are only interested in
$\phi^{-1}(0,0,-1) \cap B$, which is precisely the set of anticomplex
points in the original disk $D$. Each such anticomplex point
$(x_i,y_i)$ becomes a circle $L_i = \{(x_i,y_i,z,t) | \sqrt{z^2+t^2} =
\e \}$ in $E$ with framing $0$, oriented against the orientation of
$S^1_\e = \p D^2_\e$ if $(x_i,y_i)$ is elliptic, and with the
orientation of $S^1_\e$ if $(x_i,y_i)$ is hyperbolic. Thus the
complete oriented framed link $\psi^{-1}(0,0,-1)$ is a $-1$--framed
unknot $L_0$ with $n$ $0$--framed meridians $L_1, \ldots, L_n$ (where
$D$ has $n$ anticomplex points), with $\lk(L_0,L_i) = +1$ (respctively $-1$)
if $(x_i,y_i)$ is elliptic (respectively hyperbolic); see
Figure~\ref{F:Pontrjagin}. This is framed cobordant to an unknot with
framing equal to the sum of the entries in the linking matrix, which
is $-1 + 2e_- - 2h_-$.
\end{proof}

\section{Cocycles that guide our construction}
\label{cocycles}

Recall that $X_1$ is the union of the $0$-- and $1$--handles of $X$, and
that we have a triple $(\o_1, J_1, \xi_1)$ on $(X_1,Y_1)$, with zero
circles $Z = Z_1 \cup \ldots \cup Z_n \subset \interior(X_1)$. As in
the introduction, fix a trivialization $\tau$ of $\xi_1$; this can be
done because $c_1(\xi_1) = c_1(J_1|Y_1) = c_1(J_1)|Y_1 = 0$, since
$H^2(X_1;\Z) = 0$.

Now suppose that $K_1, \ldots, K_q$, for some $q$, are disjoint
Legendrian knots in $(Y_1,\xi_1)$, with rotation numbers $\rot(K_i)$
measured relative to $\tau$. Attach symplectic $2$--handles $H^2_1,
\ldots, H^2_q$ along $K_1, \ldots, K_q$ to produce a $4$--manifold $W$
with boundary and a triple $(\o_W, J_W, \xi_W)$ on $(W, \p W)$ that
extends $(\o_1, J_1, \xi_1)$. $W$ deformation retracts onto a
$2$--complex, with $2$--cells $C^2_1, \ldots, C^2_q$ corresponding to
the $2$--handles $H^2_1, \ldots, H^2_q$, with $1$--cells coming from the
$1$--handles of $X_1$, and with one $0$--cell at the center of the
$0$--handle. We may assume that $Z$ misses the $2$--skeleton. Thus we
get a $2$--cochain $r$ on the $2$--skeleton given by $r(C^2_i) =
\rot(K_i)$. This is trivially a cocycle on the $2$--skeleton (there are
no $3$--cells), and therefore defines a cohomology class $[r] \in
H^2(W;\Z)$. The following result is a slight generalization of a
standard fact (see~\cite{gompf}) relating rotation numbers and Chern
classes.

\begin{lemma}
\label{cocyclec1}
$[r] = c_1(J_W)$
\end{lemma}

\begin{proof}
In fact we will show that $r$ is a cocycle given by trivializing $J_W$
(by which we mean the $\C^2$--bundle defined by $J_W$) over the
$1$--skeleton and measuring the obstruction to extending this
trivialization over each $2$--cell. $\tau$ together with the outward
normal $\nu$ to $\p X_1$ defines a trivialization $(\tau,\nu)$ of
$J_W$ over $\p X_1$. By an isotopy of $W$ we can arrange that the
$1$--skeleton lies in $\p X_1$ and that the $2$--skeleton misses the
interior of $W_1$. Thus we take $(\tau,\nu)$ as the starting
trivialization of $J_W$ over the $1$--skeleton. Clearly $(\tau,\nu)$ extends
across as much of each $2$--cell $C^2_i$ as lies in $\p W_1$; the rest
of $C^2_i$ is simply the core of $H^2_i$ and our task is to show that
the obstruction to extending $(\tau,\nu)$ across $H^2_i$ is exactly
$\rot(K_i)$. Let $\kappa$ be the unit tangent vector to $K_i$; in the
proof of Proposition~2.3 in~\cite{gompf} it is shown that
$(\kappa,\nu)$ extends across $H_i$. From this we see that the winding
of $\kappa$ with respect to $\tau$ is precisely the obstruction in
question.
\end{proof}

In our main construction outlined in the introduction, we use this
lemma three times, once to recognize $c_1(J'_2)$ in our ``false
start'' construction of $(\o'_2,J'_2,\xi'_2)$, and then to
construct the correct $(\o_N,J_N,\xi_N)$ and its extension to
$(\o_2,J_2,\xi_2)$. Recall that the cocycle $x$ used to
determine the rotation numbers which produce
$(\o_2,J_2,\xi_2)$ is an honest cocycle for the full
$4$--complex decomposition of $X$.

\begin{lemma}
\label{cocycleextends}
Because $x$ is a cocycle, $J_2$ extends over the $3$--skeleton
$X$.
\end{lemma}

\begin{proof}
Consider a $3$--cell $C^3_j$. Since $x$ is a cocycle, $x(\p C^3_j) =
0$, which, by the same argument as in the preceding lemma, means that
$c_1(J_2|\p C^3_j) = 0 \in H^2(\p C^3_j;\Z)$. Trivialize $TX|C^3_j$ as
$\R^4 \times C^3_j$, with standard basis vectors $e_1, e_2, e_3, e_4$
for $\R^4$; almost complex structures on $C^3_j$ and on $\p C^3_j$ are
then determined up to homotopy by the map $J(e_4)\co  \p C^3_j
\rightarrow S^2$, where $S^2$ is the unit $(e_1,e_2,e_3)$--sphere. For
any almost complex structure $J$ on $\p C^3_j$, it is not hard to see
that $c_1(J) \in H^2(\p C^3_j;\Z) = \Z$ is twice the degree of
$J(e_4)$. (This is most easily seen by noting that the $C^2$ bundle
over $\p C^3_j$ is the Whitney sum of a trivial complex line bundle
spanned by $e_4$ and a complex line bundle $\xi$; but $\xi$ is visibly
the tangent bundle of $S^2$ when $J(e_4)$ is the identity map, so
$c_1(J)$ is twice the degree in the generating case.) Thus in our case
$J_2(e_4)$ has degree $0$ and therefore extends over $C^3_j$.
\end{proof}

\section{Plane fields and $\text{spin}^\C$ structures}
\label{planes}

Recall that in our main construction we had constructed two almost
complex structure $J_E$ and $J_*|E$ on $E$, with $c_1(J_E) = c_1(J_*|E)$.

\begin{lemma}
\label{Jhomotopy}
This implies that $J_E$ is homotopic to $J_*|E$.
\end{lemma}

\begin{proof}
Pick a nowhere zero section $\tau$ of $TE$ (an orientable
$\R^4$--bundle over a $2$--complex always has a nonzero section) and let
$\mathcal{J}$ be the unit $S^2$ bundle orthogonal to $\tau$; up to
homotopy we can view $J_E$ and $J_*$ as sections of
$\mathcal{J}$. (This is simply the bundle of local almost complex
structures.) We can homotope $J_E$ to agree with $J_*$ over the
$1$--skeleton. Over each $2$--cell $C^2_i$, we can trivialize
$\mathcal{J}$ as $S^2 \times C^2_i$ so that $J_E$ is a constant
section. Then $J_*$ becomes a map from $C^2_i$ to $S^2$ which is
constant on $\p C^2_i$, which we can think of as a map from $S^2$ to
$S^2$, and read off the degree of this map. This gives a cocycle and,
much as in the preceding section, it is not hard to see that twice
this cocycle represents $c_1(J_E) - c_1(J_*)$. Thus, because
$H^2(E;\Z)$ has no $2$--torsion and $c_1(J_E) - c_1(J_*) = 0$, we see
that this cocycle is a coboundary, which implies that we can change
our choice of trivialization over the $1$--skeleton to make the cocycle
$0$. Thus $J_E$ is homotopic to $J_*|E$.
\end{proof}

Next we show how to hit the right $\text{spin}^\C$ structure in our
construction.  Here we recall a special property of our cell
decomposition of $X$, that $\p C^2_{p+1} = C^1_1$, where $C^2_{p+1}$
is the $2$--cell associated to the $2$--handle $H^2_{p+1}$ coming from
$E$, and $C^1_1$ is a $1$--cell. Let $b$ be the $1$--cochain that is $1$
on $C^1_1$ and $0$ on everything else. Then every class $a \in
H^2(X;\Z)$ can be represented by a cocycle $z$ with $z(C^2_{p+1}) =
0$; if $z(C^2_{p+1}) \neq 0$, then replace $z$ with $z -
(z(C^2_{p+1}))\delta b$.

For our purposes we will think of a $\text{spin}^\C$ structure on $X$
as a homotopy class of almost complex structures over the $2$--skeleton
which extends over the $3$--skeleton. As mentioned above, we should
think of almost complex structures as sections of an $S^2$--bundle
$\mathcal{J}$. For our fixed class $c \in H^2(X;\Z)$, let
$\mathcal{S}_c$ be the set of all $\text{spin}^\C$ structures $s$ with
$c_1(s) = c$. It is well known that $H^2(X;\Z)$ acts freely and
transitively on the set of all $\text{spin}^\C$ structures on $X$,
changing $c_1$ by twice the cohomology class that is acting. Thus the
difference between two $\text{spin}^\C$ structures $s_0, s_1 \in
\mathcal{S}_c$ is a cohomology class of order $2$.

\begin{proposition}
\label{spinCcocycle}
Suppose that, in our construction, we used a cocycle representative
$x_0$ of $c$, and that this produced $s_0 \in
\mathcal{S}_c$. Let $a \in H^2(X;\Z)$ be the class of order $2$ which
acts on $s_0$ to give $s_1=s$. Choose a representative $z$ for $a$
with $z(C^2_{p+1}) = 0$. Then, if we repeat our construction with the
cocycle $x_1 = x_0 - 2z$ instead of $x_0$, we will produce
$s_1$ instead of $s_0$.
\end{proposition}

\begin{proof}
  Let $J_0, J_1$ be the almost complex structures produced by
  $x_1,x_0$, respectively. Recall that our construction depended on a
  fixed trivialization over the $1$--skeleton, and that we can homotope
  $J_0$ to $J_1$ over the $1$--skeleton. Again, over each $2$--cell
  $C^2_i$ we can trivialize $\mathcal{J}$ so that $J_0$ is a constant
  section. Then the action of $a$ is given precisely by changing this
  section to a section of degree $z(C^2_i)$. The difference between
  the obstructions to extending the trivialization of $J_1 = J_0$ on
  the $1$--skeleton to a trivialization of $J_1$ or $J_0$ over $C^2_i$
  is twice this degree, and thus we see that changing the rotation
  numbers by $2z$ implements this action of $a$.

It is important that $z(C^2_{p+1}) = 0$ because in our final
construction we discard $J_2$ and only use $J_N$. Thus, if we had
changed our construction on $H^2_{p+1}$, the change would not survive
to the end of the construction.
\end{proof}

\end{document}